\documentclass{gtart}

\def\ifplaintex{\expandafter\ifx\csname documentclass\endcsname\relax}

\ifplaintex 
\else
\fi

\expandafter\ifx\csname beginpicture\endcsname\relax
\expandafter\ifx\csname documentclass\endcsname\relax
\input pictex \else
\input prepictex \input pictex \input postpictex \fi\fi

\def\gt{{\mathsurround=0pt\it $\cal G\mskip-2mu$eometry \&\ 
$\cal T\!\!$opology}}        

\def\gtp{{\mathsurround=0pt\it $\cal G\mskip-2mu$eometry \&\ 
$\cal T\!\!$opology $\cal P\!$ublications}}  

\def\lognumber#1{\def\thelognumber{#1}}
\def\volumenumber#1{\def\thevolumenumber{#1}}
\def\papernumber#1{\def\thepapernumber{#1}}
\def\volumeyear#1{\def\thevolumeyear{#1}}

\def\pagenumbers#1#2{\def\startpage{#1}\def\finishpage{#2}}
\def\published#1{\def\publishdate{#1}}
\def\proposed#1{\def\theproposer{#1}}
\def\seconded#1{\def\theseconders{#1}}
\def\received#1{\def\receiveddate{#1}}
\def\revised#1{\def\reviseddate{#1}}
\def\accepted#1{\def\accepteddate{#1}}

\long\def\asciiabstract#1{\long\def\theasciiabstract{#1}}
\def\asciikeywords#1{\def\theasciikeywords{#1}}

\def\shorttitle#1{\def\theshorttitle{#1}}

\let\\\par\let\thelognumber\relax
\let\thevolumenumber\relax\let\thepapernumber\relax
\let\thevolumeyear\relax\let\thesamplenumber\relax\let\startpage\relax
\let\finishpage\relax\let\publishdate\relax\let\receiveddate\relax
\let\reviseddate\relax\let\accepteddate\relax\let\theasciititle\relax
\let\theasciiauthors\relax
\let\theasciiabstract\relax\let\theasciikeywords\relax
\let\theasciiemail\relax\let\theshortauthors\relax\let\theshorttitle\relax

\long\def\maketitlep{   

\count0=\startpage

\gt\hfill      
\beginpicture
\setcoordinatesystem units <0.33truein, 0.33truein> point at 2.2 0.9
\setplotsymbol ({$\cal G$})
\plotsymbolspacing=9truept
\circulararc 315 degrees from 0 1 center at 0 0
\setplotsymbol ({$\cal T$})
\circulararc 315 degrees from 1 -1 center at 1 0
\endpicture
\break
{\small\ifx\thesamplenumber\relax 
Volume \else Sample
\fi\thevolumenumber\ (\thevolumeyear)
\startpage--\finishpage\nl
Published: \publishdate}
\vglue 0.5truein plus 0.4fil minus 0.1truein

{\parskip=0pt\leftskip 0pt plus 1fil\def\\{\par\smallskip}{\ifplaintex\large
\else\Large\fi\bf\thetitle}\par\medskip}   

\vglue 0pt plus 0.1fil 

{\parskip=0pt\leftskip 0pt plus 1fil\def\\{\par}{\sc\theauthors}
\par\medskip}

\vglue 0pt plus 0.1fil 

{\small\parskip=0pt\let\newline\\
{\leftskip 0pt plus 1fil\def\\{\par}{\sl\theaddress}\par}
\expandafter\ifx\theemail\relax    
\relax\else\vglue 5pt plus 0.02fil minus 2pt\def\\{\stdspace{\rm 
and}\stdspace} 
\cl{Email:\stdspace\tt\theemail}\fi
\ifx\theurl\relax                  
\relax\else\vglue 5pt plus 0.02fil minus 2pt\def\\{\stdspace{\rm 
and}\stdspace}
\cl{URL:\stdspace\tt\theurl}\fi\par}

\vglue 7pt plus 0.3fil minus 3pt

{\bf Abstract}
\vglue 5pt plus 0.1fil minus 2pt

\theabstract

\vglue 7pt plus 0.3fil minus 3pt

{\bf AMS Classification numbers}\quad Primary:\quad \theprimaryclass

Secondary:\quad \thesecondaryclass

\vglue 5pt plus 0.3fil minus 2pt

{\bf Keywords}\quad \thekeywords

\vglue 10pt plus 0.5fil minus 5pt

{\small  Proposed: \theproposer\hfill Received: \receiveddate\nl
Seconded: \theseconders\hfill 
\ifx\reviseddate\relax                         
Accepted: \accepteddate                        
\else
Revised: \reviseddate                          
\fi}
\eject
}       

\let\maketitlepage\maketitlep
\let\maketitle\maketitlepage

\font\phead=cmsl9 scaled 950
\font=cmsl9 scaled 1050
\font\pnum=cmbx10 scaled 913
\font\lnum=cmbx10 
\font\pfoot=cmsl9 scaled 950
\font=cmsl9 scaled 1050
\ifplaintex
\headline{\vbox to 0pt{\vskip -4.5mm\line{\small\phead\ifnum
\count0=\startpage ISSN 1364-0380 (on line)
1465-3060 (printed) \hfill {\pnum\folio}\else\ifodd\count0\def\\{ }%
\ifx\theshorttitle\relax\thetitle\else\theshorttitle\fi\hfill{\pnum\folio}
\else\def\\{ and }{\pnum\folio}\hfill\ifx\theshortauthors\relax\theauthors
\else\theshortauthors\fi\fi\fi}\vss}}
\footline{\vbox to 0pt{\vglue 0mm\line{\small\pfoot\ifnum\count0=\startpage
\copyright\ \gtp\hfill\else
\gt, Volume \thevolumenumber\ (\thevolumeyear)\hfill\fi}\vss
}}
\else
\makeatletter
\def\@oddhead{{\small\lhead\ifnum\count0=\startpage ISSN 1364-0380 (on line)
1465-3060 (printed) \hfill {\lnum\number\count0}\else\ifodd\count0
\def\\{ }\ifx\theshorttitle\relax \thetitle \else\theshorttitle\fi\hfill
{\lnum\number\count0}\else\def\\{ and }{\lnum\number\count0}
\hfill\ifx\theshortauthors\relax 
\theauthors\else\theshortauthors\fi\fi\fi}}\def\@evenhead{\@oddhead}
\def\@oddfoot{\small\lfoot\ifnum\count0=\startpage\copyright\ \gtp\hfill\else
\gt, Volume \thevolumenumber\ (\thevolumeyear)\hfill\fi}
\def\@evenfoot{\@oddfoot}
\makeatother
\fi

\newwrite\gtoutfile
\long\gdef\makeheadfile{  
{\def\\{, }\def\s{ }
\immediate\openout\gtoutfile head.xxx
\immediate\write\gtoutfile{To: math@arxiv.org}
\immediate\write\gtoutfile{Subject: put or rep NNNNN:pppp}
\immediate\write\gtoutfile{--text follows this line--}
\immediate\write\gtoutfile{Proxy-for: \ifx\theasciiauthors\relax
\theauthors\else\theasciiauthors\fi\s<\ifx\theasciiemail\relax\theemail\else\theasciiemail\fi>}
\immediate\write\gtoutfile{\noexpand\\}
\immediate\write\gtoutfile{Authors: \ifx\theasciiauthors\relax
\theauthors\else\theasciiauthors\fi}
\immediate\write\gtoutfile{Title: \ifx\theasciititle\relax
\thetitle\else\theasciititle\fi}
\immediate\write\gtoutfile{Subj-class: GT or SG or MG etc}
\immediate\write\gtoutfile{MSC-class: \theprimaryclass\ifx\thesecondaryclass\relax\else, \thesecondaryclass\fi}
\immediate\write\gtoutfile{Journal-ref: Geom. Topol. \thevolumenumber
(\thevolumeyear) \startpage-\finishpage}
\immediate\write\gtoutfile{Comments: Published by Geometry and Topology at}
\immediate\write\gtoutfile{\s\s http://www.maths.warwick.ac.uk/gt/GTVol\thevolumenumber/paper\thepapernumber.abs.html}
\immediate\write\gtoutfile{\noexpand\\}
\immediate\write\gtoutfile{}
\ifx\theasciiabstract\relax
\immediate\write\gtoutfile{\theabstract}\else
\immediate\write\gtoutfile{\theasciiabstract}\fi
\immediate\write\gtoutfile{}
\immediate\write\gtoutfile{\noexpand\\}
\immediate\write\gtoutfile{}
\immediate\closeout\gtoutfile}}  

\def\maketitlepage{\maketitlep\makeheadfile}
\let\maketitle\maketitlepage

\def\ifplaintex{\expandafter\ifx\csname documentclass\endcsname\relax}

\ifplaintex 
\else
\fi

\expandafter\ifx\csname beginpicture\endcsname\relax
\expandafter\ifx\csname documentclass\endcsname\relax
\input pictex \else
\input prepictex \input pictex \input postpictex \fi\fi

\def\gt{{\mathsurround=0pt\it $\cal G\mskip-2mu$eometry \&\ 
$\cal T\!\!$opology}}        

\def\gtp{{\mathsurround=0pt\it $\cal G\mskip-2mu$eometry \&\ 
$\cal T\!\!$opology $\cal P\!$ublications}}  

\def\lognumber#1{\def\thelognumber{#1}}
\def\volumenumber#1{\def\thevolumenumber{#1}}
\def\papernumber#1{\def\thepapernumber{#1}}
\def\volumeyear#1{\def\thevolumeyear{#1}}

\def\pagenumbers#1#2{\def\startpage{#1}\def\finishpage{#2}}
\def\published#1{\def\publishdate{#1}}
\def\proposed#1{\def\theproposer{#1}}
\def\seconded#1{\def\theseconders{#1}}
\def\received#1{\def\receiveddate{#1}}
\def\revised#1{\def\reviseddate{#1}}
\def\accepted#1{\def\accepteddate{#1}}

\long\def\asciiabstract#1{\long\def\theasciiabstract{#1}}
\def\asciikeywords#1{\def\theasciikeywords{#1}}

\def\shorttitle#1{\def\theshorttitle{#1}}

\let\\\par\let\thelognumber\relax
\let\thevolumenumber\relax\let\thepapernumber\relax
\let\thevolumeyear\relax\let\thesamplenumber\relax\let\startpage\relax
\let\finishpage\relax\let\publishdate\relax\let\receiveddate\relax
\let\reviseddate\relax\let\accepteddate\relax\let\theasciititle\relax
\let\theasciiauthors\relax
\let\theasciiabstract\relax\let\theasciikeywords\relax
\let\theasciiemail\relax\let\theshortauthors\relax\let\theshorttitle\relax

\long\def\maketitlep{   

\count0=\startpage

\gt\hfill      
\beginpicture
\setcoordinatesystem units <0.33truein, 0.33truein> point at 2.2 0.9
\setplotsymbol ({$\cal G$})
\plotsymbolspacing=9truept
\circulararc 315 degrees from 0 1 center at 0 0
\setplotsymbol ({$\cal T$})
\circulararc 315 degrees from 1 -1 center at 1 0
\endpicture
\break
{\small\ifx\thesamplenumber\relax 
Volume \else Sample
\fi\thevolumenumber\ (\thevolumeyear)
\startpage--\finishpage\nl
Published: \publishdate}
\vglue 0.5truein plus 0.4fil minus 0.1truein

{\parskip=0pt\leftskip 0pt plus 1fil\def\\{\par\smallskip}{\ifplaintex\large
\else\Large\fi\bf\thetitle}\par\medskip}   

\vglue 0pt plus 0.1fil 

{\parskip=0pt\leftskip 0pt plus 1fil\def\\{\par}{\sc\theauthors}
\par\medskip}

\vglue 0pt plus 0.1fil 

{\small\parskip=0pt\let\newline\\
{\leftskip 0pt plus 1fil\def\\{\par}{\sl\theaddress}\par}
\expandafter\ifx\theemail\relax    
\relax\else\vglue 5pt plus 0.02fil minus 2pt\def\\{\stdspace{\rm 
and}\stdspace} 
\cl{Email:\stdspace\tt\theemail}\fi
\ifx\theurl\relax                  
\relax\else\vglue 5pt plus 0.02fil minus 2pt\def\\{\stdspace{\rm 
and}\stdspace}
\cl{URL:\stdspace\tt\theurl}\fi\par}

\vglue 7pt plus 0.3fil minus 3pt

{\bf Abstract}
\vglue 5pt plus 0.1fil minus 2pt

\theabstract

\vglue 7pt plus 0.3fil minus 3pt

{\bf AMS Classification numbers}\quad Primary:\quad \theprimaryclass

Secondary:\quad \thesecondaryclass

\vglue 5pt plus 0.3fil minus 2pt

{\bf Keywords}\quad \thekeywords

\vglue 10pt plus 0.5fil minus 5pt

{\small  Proposed: \theproposer\hfill Received: \receiveddate\nl
Seconded: \theseconders\hfill 
\ifx\reviseddate\relax                         
Accepted: \accepteddate                        
\else
Revised: \reviseddate                          
\fi}
\eject
}       

\let\maketitlepage\maketitlep
\let\maketitle\maketitlepage

\font\phead=cmsl9 scaled 950
\font=cmsl9 scaled 1050
\font\pnum=cmbx10 scaled 913
\font\lnum=cmbx10 
\font\pfoot=cmsl9 scaled 950
\font=cmsl9 scaled 1050
\ifplaintex
\headline{\vbox to 0pt{\vskip -4.5mm\line{\small\phead\ifnum
\count0=\startpage ISSN 1364-0380 (on line)
1465-3060 (printed) \hfill {\pnum\folio}\else\ifodd\count0\def\\{ }%
\ifx\theshorttitle\relax\thetitle\else\theshorttitle\fi\hfill{\pnum\folio}
\else\def\\{ and }{\pnum\folio}\hfill\ifx\theshortauthors\relax\theauthors
\else\theshortauthors\fi\fi\fi}\vss}}
\footline{\vbox to 0pt{\vglue 0mm\line{\small\pfoot\ifnum\count0=\startpage
\copyright\ \gtp\hfill\else
\gt, Volume \thevolumenumber\ (\thevolumeyear)\hfill\fi}\vss
}}
\else
\makeatletter
\def\@oddhead{{\small\lhead\ifnum\count0=\startpage ISSN 1364-0380 (on line)
1465-3060 (printed) \hfill {\lnum\number\count0}\else\ifodd\count0
\def\\{ }\ifx\theshorttitle\relax \thetitle \else\theshorttitle\fi\hfill
{\lnum\number\count0}\else\def\\{ and }{\lnum\number\count0}
\hfill\ifx\theshortauthors\relax 
\theauthors\else\theshortauthors\fi\fi\fi}}\def\@evenhead{\@oddhead}
\def\@oddfoot{\small\lfoot\ifnum\count0=\startpage\copyright\ \gtp\hfill\else
\gt, Volume \thevolumenumber\ (\thevolumeyear)\hfill\fi}
\def\@evenfoot{\@oddfoot}
\makeatother
\fi

\newwrite\gtoutfile
\long\gdef\makeheadfile{  
{\def\\{, }\def\s{ }
\immediate\openout\gtoutfile head.xxx
\immediate\write\gtoutfile{To: math@arxiv.org}
\immediate\write\gtoutfile{Subject: put or rep NNNNN:pppp}
\immediate\write\gtoutfile{--text follows this line--}
\immediate\write\gtoutfile{Proxy-for: \ifx\theasciiauthors\relax
\theauthors\else\theasciiauthors\fi\s<\ifx\theasciiemail\relax\theemail\else\theasciiemail\fi>}
\immediate\write\gtoutfile{\noexpand\\}
\immediate\write\gtoutfile{Authors: \ifx\theasciiauthors\relax
\theauthors\else\theasciiauthors\fi}
\immediate\write\gtoutfile{Title: \ifx\theasciititle\relax
\thetitle\else\theasciititle\fi}
\immediate\write\gtoutfile{Subj-class: GT or SG or MG etc}
\immediate\write\gtoutfile{MSC-class: \theprimaryclass\ifx\thesecondaryclass\relax\else, \thesecondaryclass\fi}
\immediate\write\gtoutfile{Journal-ref: Geom. Topol. \thevolumenumber
(\thevolumeyear) \startpage-\finishpage}
\immediate\write\gtoutfile{Comments: Published by Geometry and Topology at}
\immediate\write\gtoutfile{\s\s http://www.maths.warwick.ac.uk/gt/GTVol\thevolumenumber/paper\thepapernumber.abs.html}
\immediate\write\gtoutfile{\noexpand\\}
\immediate\write\gtoutfile{}
\ifx\theasciiabstract\relax
\immediate\write\gtoutfile{\theabstract}\else
\immediate\write\gtoutfile{\theasciiabstract}\fi
\immediate\write\gtoutfile{}
\immediate\write\gtoutfile{\noexpand\\}
\immediate\write\gtoutfile{}
\immediate\closeout\gtoutfile}}  

\def\maketitlepage{\maketitlep\makeheadfile}
\let\maketitle\maketitlepage

\def\ifplaintex{\expandafter\ifx\csname documentclass\endcsname\relax}

\ifplaintex 
\else
\fi

\expandafter\ifx\csname beginpicture\endcsname\relax
\expandafter\ifx\csname documentclass\endcsname\relax
\input pictex \else
\input prepictex \input pictex \input postpictex \fi\fi

\def\gt{{\mathsurround=0pt\it $\cal G\mskip-2mu$eometry \&\ 
$\cal T\!\!$opology}}        

\def\gtp{{\mathsurround=0pt\it $\cal G\mskip-2mu$eometry \&\ 
$\cal T\!\!$opology $\cal P\!$ublications}}  

\def\lognumber#1{\def\thelognumber{#1}}
\def\volumenumber#1{\def\thevolumenumber{#1}}
\def\papernumber#1{\def\thepapernumber{#1}}
\def\volumeyear#1{\def\thevolumeyear{#1}}

\def\pagenumbers#1#2{\def\startpage{#1}\def\finishpage{#2}}
\def\published#1{\def\publishdate{#1}}
\def\proposed#1{\def\theproposer{#1}}
\def\seconded#1{\def\theseconders{#1}}
\def\received#1{\def\receiveddate{#1}}
\def\revised#1{\def\reviseddate{#1}}
\def\accepted#1{\def\accepteddate{#1}}

\long\def\asciiabstract#1{\long\def\theasciiabstract{#1}}
\def\asciikeywords#1{\def\theasciikeywords{#1}}

\def\shorttitle#1{\def\theshorttitle{#1}}

\let\\\par\let\thelognumber\relax
\let\thevolumenumber\relax\let\thepapernumber\relax
\let\thevolumeyear\relax\let\thesamplenumber\relax\let\startpage\relax
\let\finishpage\relax\let\publishdate\relax\let\receiveddate\relax
\let\reviseddate\relax\let\accepteddate\relax\let\theasciititle\relax
\let\theasciiauthors\relax
\let\theasciiabstract\relax\let\theasciikeywords\relax
\let\theasciiemail\relax\let\theshortauthors\relax\let\theshorttitle\relax

\long\def\maketitlep{   

\count0=\startpage

\gt\hfill      
\beginpicture
\setcoordinatesystem units <0.33truein, 0.33truein> point at 2.2 0.9
\setplotsymbol ({$\cal G$})
\plotsymbolspacing=9truept
\circulararc 315 degrees from 0 1 center at 0 0
\setplotsymbol ({$\cal T$})
\circulararc 315 degrees from 1 -1 center at 1 0
\endpicture
\break
{\small\ifx\thesamplenumber\relax 
Volume \else Sample
\fi\thevolumenumber\ (\thevolumeyear)
\startpage--\finishpage\nl
Published: \publishdate}
\vglue 0.5truein plus 0.4fil minus 0.1truein

{\parskip=0pt\leftskip 0pt plus 1fil\def\\{\par\smallskip}{\ifplaintex\large
\else\Large\fi\bf\thetitle}\par\medskip}   

\vglue 0pt plus 0.1fil 

{\parskip=0pt\leftskip 0pt plus 1fil\def\\{\par}{\sc\theauthors}
\par\medskip}

\vglue 0pt plus 0.1fil 

{\small\parskip=0pt\let\newline\\
{\leftskip 0pt plus 1fil\def\\{\par}{\sl\theaddress}\par}
\expandafter\ifx\theemail\relax    
\relax\else\vglue 5pt plus 0.02fil minus 2pt\def\\{\stdspace{\rm 
and}\stdspace} 
\cl{Email:\stdspace\tt\theemail}\fi
\ifx\theurl\relax                  
\relax\else\vglue 5pt plus 0.02fil minus 2pt\def\\{\stdspace{\rm 
and}\stdspace}
\cl{URL:\stdspace\tt\theurl}\fi\par}

\vglue 7pt plus 0.3fil minus 3pt

{\bf Abstract}
\vglue 5pt plus 0.1fil minus 2pt

\theabstract

\vglue 7pt plus 0.3fil minus 3pt

{\bf AMS Classification numbers}\quad Primary:\quad \theprimaryclass

Secondary:\quad \thesecondaryclass

\vglue 5pt plus 0.3fil minus 2pt

{\bf Keywords}\quad \thekeywords

\vglue 10pt plus 0.5fil minus 5pt

{\small  Proposed: \theproposer\hfill Received: \receiveddate\nl
Seconded: \theseconders\hfill 
\ifx\reviseddate\relax                         
Accepted: \accepteddate                        
\else
Revised: \reviseddate                          
\fi}
\eject
}       

\let\maketitlepage\maketitlep
\let\maketitle\maketitlepage

\font\phead=cmsl9 scaled 950
\font=cmsl9 scaled 1050
\font\pnum=cmbx10 scaled 913
\font\lnum=cmbx10 
\font\pfoot=cmsl9 scaled 950
\font=cmsl9 scaled 1050
\ifplaintex
\headline{\vbox to 0pt{\vskip -4.5mm\line{\small\phead\ifnum
\count0=\startpage ISSN 1364-0380 (on line)
1465-3060 (printed) \hfill {\pnum\folio}\else\ifodd\count0\def\\{ }%
\ifx\theshorttitle\relax\thetitle\else\theshorttitle\fi\hfill{\pnum\folio}
\else\def\\{ and }{\pnum\folio}\hfill\ifx\theshortauthors\relax\theauthors
\else\theshortauthors\fi\fi\fi}\vss}}
\footline{\vbox to 0pt{\vglue 0mm\line{\small\pfoot\ifnum\count0=\startpage
\copyright\ \gtp\hfill\else
\gt, Volume \thevolumenumber\ (\thevolumeyear)\hfill\fi}\vss
}}
\else
\makeatletter
\def\@oddhead{{\small\lhead\ifnum\count0=\startpage ISSN 1364-0380 (on line)
1465-3060 (printed) \hfill {\lnum\number\count0}\else\ifodd\count0
\def\\{ }\ifx\theshorttitle\relax \thetitle \else\theshorttitle\fi\hfill
{\lnum\number\count0}\else\def\\{ and }{\lnum\number\count0}
\hfill\ifx\theshortauthors\relax 
\theauthors\else\theshortauthors\fi\fi\fi}}\def\@evenhead{\@oddhead}
\def\@oddfoot{\small\lfoot\ifnum\count0=\startpage\copyright\ \gtp\hfill\else
\gt, Volume \thevolumenumber\ (\thevolumeyear)\hfill\fi}
\def\@evenfoot{\@oddfoot}
\makeatother
\fi

\newwrite\gtoutfile
\long\gdef\makeheadfile{  
{\def\\{, }\def\s{ }
\immediate\openout\gtoutfile head.xxx
\immediate\write\gtoutfile{To: math@arxiv.org}
\immediate\write\gtoutfile{Subject: put or rep NNNNN:pppp}
\immediate\write\gtoutfile{--text follows this line--}
\immediate\write\gtoutfile{Proxy-for: \ifx\theasciiauthors\relax
\theauthors\else\theasciiauthors\fi\s<\ifx\theasciiemail\relax\theemail\else\theasciiemail\fi>}
\immediate\write\gtoutfile{\noexpand\\}
\immediate\write\gtoutfile{Authors: \ifx\theasciiauthors\relax
\theauthors\else\theasciiauthors\fi}
\immediate\write\gtoutfile{Title: \ifx\theasciititle\relax
\thetitle\else\theasciititle\fi}
\immediate\write\gtoutfile{Subj-class: GT or SG or MG etc}
\immediate\write\gtoutfile{MSC-class: \theprimaryclass\ifx\thesecondaryclass\relax\else, \thesecondaryclass\fi}
\immediate\write\gtoutfile{Journal-ref: Geom. Topol. \thevolumenumber
(\thevolumeyear) \startpage-\finishpage}
\immediate\write\gtoutfile{Comments: Published by Geometry and Topology at}
\immediate\write\gtoutfile{\s\s http://www.maths.warwick.ac.uk/gt/GTVol\thevolumenumber/paper\thepapernumber.abs.html}
\immediate\write\gtoutfile{\noexpand\\}
\immediate\write\gtoutfile{}
\ifx\theasciiabstract\relax
\immediate\write\gtoutfile{\theabstract}\else
\immediate\write\gtoutfile{\theasciiabstract}\fi
\immediate\write\gtoutfile{}
\immediate\write\gtoutfile{\noexpand\\}
\immediate\write\gtoutfile{}
\immediate\closeout\gtoutfile}}  

\def\maketitlepage{\maketitlep\makeheadfile}
\let\maketitle\maketitlepage

\def\ifplaintex{\expandafter\ifx\csname documentclass\endcsname\relax}

\ifplaintex 
\else
\fi

\expandafter\ifx\csname beginpicture\endcsname\relax
\expandafter\ifx\csname documentclass\endcsname\relax
\input pictex \else
\input prepictex \input pictex \input postpictex \fi\fi

\def\gt{{\mathsurround=0pt\it $\cal G\mskip-2mu$eometry \&\ 
$\cal T\!\!$opology}}        

\def\gtp{{\mathsurround=0pt\it $\cal G\mskip-2mu$eometry \&\ 
$\cal T\!\!$opology $\cal P\!$ublications}}  

\def\lognumber#1{\def\thelognumber{#1}}
\def\volumenumber#1{\def\thevolumenumber{#1}}
\def\papernumber#1{\def\thepapernumber{#1}}
\def\volumeyear#1{\def\thevolumeyear{#1}}

\def\pagenumbers#1#2{\def\startpage{#1}\def\finishpage{#2}}
\def\published#1{\def\publishdate{#1}}
\def\proposed#1{\def\theproposer{#1}}
\def\seconded#1{\def\theseconders{#1}}
\def\received#1{\def\receiveddate{#1}}
\def\revised#1{\def\reviseddate{#1}}
\def\accepted#1{\def\accepteddate{#1}}

\long\def\asciiabstract#1{\long\def\theasciiabstract{#1}}
\def\asciikeywords#1{\def\theasciikeywords{#1}}

\def\shorttitle#1{\def\theshorttitle{#1}}

\let\\\par\let\thelognumber\relax
\let\thevolumenumber\relax\let\thepapernumber\relax
\let\thevolumeyear\relax\let\thesamplenumber\relax\let\startpage\relax
\let\finishpage\relax\let\publishdate\relax\let\receiveddate\relax
\let\reviseddate\relax\let\accepteddate\relax\let\theasciititle\relax
\let\theasciiauthors\relax
\let\theasciiabstract\relax\let\theasciikeywords\relax
\let\theasciiemail\relax\let\theshortauthors\relax\let\theshorttitle\relax

\long\def\maketitlep{   

\count0=\startpage

\gt\hfill      
\beginpicture
\setcoordinatesystem units <0.33truein, 0.33truein> point at 2.2 0.9
\setplotsymbol ({$\cal G$})
\plotsymbolspacing=9truept
\circulararc 315 degrees from 0 1 center at 0 0
\setplotsymbol ({$\cal T$})
\circulararc 315 degrees from 1 -1 center at 1 0
\endpicture
\break
{\small\ifx\thesamplenumber\relax 
Volume \else Sample
\fi\thevolumenumber\ (\thevolumeyear)
\startpage--\finishpage\nl
Published: \publishdate}
\vglue 0.5truein plus 0.4fil minus 0.1truein

{\parskip=0pt\leftskip 0pt plus 1fil\def\\{\par\smallskip}{\ifplaintex\large
\else\Large\fi\bf\thetitle}\par\medskip}   

\vglue 0pt plus 0.1fil 

{\parskip=0pt\leftskip 0pt plus 1fil\def\\{\par}{\sc\theauthors}
\par\medskip}

\vglue 0pt plus 0.1fil 

{\small\parskip=0pt\let\newline\\
{\leftskip 0pt plus 1fil\def\\{\par}{\sl\theaddress}\par}
\expandafter\ifx\theemail\relax    
\relax\else\vglue 5pt plus 0.02fil minus 2pt\def\\{\stdspace{\rm 
and}\stdspace} 
\cl{Email:\stdspace\tt\theemail}\fi
\ifx\theurl\relax                  
\relax\else\vglue 5pt plus 0.02fil minus 2pt\def\\{\stdspace{\rm 
and}\stdspace}
\cl{URL:\stdspace\tt\theurl}\fi\par}

\vglue 7pt plus 0.3fil minus 3pt

{\bf Abstract}
\vglue 5pt plus 0.1fil minus 2pt

\theabstract

\vglue 7pt plus 0.3fil minus 3pt

{\bf AMS Classification numbers}\quad Primary:\quad \theprimaryclass

Secondary:\quad \thesecondaryclass

\vglue 5pt plus 0.3fil minus 2pt

{\bf Keywords}\quad \thekeywords

\vglue 10pt plus 0.5fil minus 5pt

{\small  Proposed: \theproposer\hfill Received: \receiveddate\nl
Seconded: \theseconders\hfill 
\ifx\reviseddate\relax                         
Accepted: \accepteddate                        
\else
Revised: \reviseddate                          
\fi}
\eject
}       

\let\maketitlepage\maketitlep
\let\maketitle\maketitlepage

\font\phead=cmsl9 scaled 950
\font=cmsl9 scaled 1050
\font\pnum=cmbx10 scaled 913
\font\lnum=cmbx10 
\font\pfoot=cmsl9 scaled 950
\font=cmsl9 scaled 1050
\ifplaintex
\headline{\vbox to 0pt{\vskip -4.5mm\line{\small\phead\ifnum
\count0=\startpage ISSN 1364-0380 (on line)
1465-3060 (printed) \hfill {\pnum\folio}\else\ifodd\count0\def\\{ }%
\ifx\theshorttitle\relax\thetitle\else\theshorttitle\fi\hfill{\pnum\folio}
\else\def\\{ and }{\pnum\folio}\hfill\ifx\theshortauthors\relax\theauthors
\else\theshortauthors\fi\fi\fi}\vss}}
\footline{\vbox to 0pt{\vglue 0mm\line{\small\pfoot\ifnum\count0=\startpage
\copyright\ \gtp\hfill\else
\gt, Volume \thevolumenumber\ (\thevolumeyear)\hfill\fi}\vss
}}
\else
\makeatletter
\def\@oddhead{{\small\lhead\ifnum\count0=\startpage ISSN 1364-0380 (on line)
1465-3060 (printed) \hfill {\lnum\number\count0}\else\ifodd\count0
\def\\{ }\ifx\theshorttitle\relax \thetitle \else\theshorttitle\fi\hfill
{\lnum\number\count0}\else\def\\{ and }{\lnum\number\count0}
\hfill\ifx\theshortauthors\relax 
\theauthors\else\theshortauthors\fi\fi\fi}}\def\@evenhead{\@oddhead}
\def\@oddfoot{\small\lfoot\ifnum\count0=\startpage\copyright\ \gtp\hfill\else
\gt, Volume \thevolumenumber\ (\thevolumeyear)\hfill\fi}
\def\@evenfoot{\@oddfoot}
\makeatother
\fi

\newwrite\gtoutfile
\long\gdef\makeheadfile{  
{\def\\{, }\def\s{ }
\immediate\openout\gtoutfile head.xxx
\immediate\write\gtoutfile{To: math@arxiv.org}
\immediate\write\gtoutfile{Subject: put or rep NNNNN:pppp}
\immediate\write\gtoutfile{--text follows this line--}
\immediate\write\gtoutfile{Proxy-for: \ifx\theasciiauthors\relax
\theauthors\else\theasciiauthors\fi\s<\ifx\theasciiemail\relax\theemail\else\theasciiemail\fi>}
\immediate\write\gtoutfile{\noexpand\\}
\immediate\write\gtoutfile{Authors: \ifx\theasciiauthors\relax
\theauthors\else\theasciiauthors\fi}
\immediate\write\gtoutfile{Title: \ifx\theasciititle\relax
\thetitle\else\theasciititle\fi}
\immediate\write\gtoutfile{Subj-class: GT or SG or MG etc}
\immediate\write\gtoutfile{MSC-class: \theprimaryclass\ifx\thesecondaryclass\relax\else, \thesecondaryclass\fi}
\immediate\write\gtoutfile{Journal-ref: Geom. Topol. \thevolumenumber
(\thevolumeyear) \startpage-\finishpage}
\immediate\write\gtoutfile{Comments: Published by Geometry and Topology at}
\immediate\write\gtoutfile{\s\s http://www.maths.warwick.ac.uk/gt/GTVol\thevolumenumber/paper\thepapernumber.abs.html}
\immediate\write\gtoutfile{\noexpand\\}
\immediate\write\gtoutfile{}
\ifx\theasciiabstract\relax
\immediate\write\gtoutfile{\theabstract}\else
\immediate\write\gtoutfile{\theasciiabstract}\fi
\immediate\write\gtoutfile{}
\immediate\write\gtoutfile{\noexpand\\}
\immediate\write\gtoutfile{}
\immediate\closeout\gtoutfile}}  

\def\maketitlepage{\maketitlep\makeheadfile}
\let\maketitle\maketitlepage

\lognumber{102}
\volumenumber{5}\papernumber{22}\volumeyear{2001}
\pagenumbers{683}{718}
\received{2 November 1999}
\revised{28 August 2001}
\accepted{26 September 2001}
\proposed{Ralph Cohen}
\seconded{Haynes Miller, Steven Ferry}
\published{26 September 2001}

\usepackage{amscd,graphicx}

\newenvironment{Relax}{\relax}{\relax}

\newtheorem{Theorem}{Theorem}[section]
\newtheorem{Lemma}[Theorem]{Lemma}

\newtheorem{Corollary}[Theorem]{Corollary}

\theoremstyle{remark}
\newtheorem{Remark}[Theorem]{Remark}
\newtheorem{Definition}[Theorem]{Definition}

\def\strut#1#2{\vrule width 0pt height #1pt depth #2pt}
\def\overset#1\to#2{{ \stackrel{#1}{#2} }}
\def\underset#1\to#2{{ {#2}_{#1} }}

\def\Cal#1{{\cal#1}}

\let\text=\mbox

\def\({ \left( }
\def\){ \right) }
\def\[{ \left[ }
\def\]{ \right] }
\def\<{ \langle }
\def\>{ \rangle }
\let\ljunk=\{
\let\rjunk=\}
\def\{{\left\ljunk}
\def\}{\right\rjunk}
\def\^{ \wedge }
\def\@{ \otimes }
\newcommand{\PSbox}[3]{\mbox{\rule{0in}{#3}\smash{\rlap{\hbox{\includegraphics{#1}}}}\hspace{#2}}}
\def\PSboxa <#1,#2> #3#4#5{\mbox{\rule{0in}{#5}\smash{\rlap{\kern #1 \raise #2\hbox{\includegraphics{#3}}}}\hspace{#4}}}
\def\R{\mbox{{\bf R}}}
\newcommand{\Z}{{\mathbf Z}}
\def\C{\mbox{{\bf C}}}
\def\H{\mbox{{\bf H}}}
\def\P{\mbox{{\bf P}}}

\def\p{\partial}

\def\X{\times}
\def\A{{\Cal A}}
\def\S{\Sigma}

\def\X{\times}

\def\ov{\overline}

\def\ind{\mbox{ind}}
\def\Ker{\mbox{ Ker }}
\def\Im{\mbox{ Im }}
\def\Hom{\mbox{Hom}}
\def\Ext{\mbox{Ext}}

\begin{document}
\title{Manifolds with singularities accepting a metric\\of positive
scalar curvature} 
\shorttitle{Manifolds with singularities of positive scalar curvature} 
\author{Boris Botvinnik} 
\address{Department of Mathematics, University of Oregon\\Eugene, 
OR 97403, USA}
\email{botvinn@math.uoregon.edu}
\begin{abstract}
	We study the question of existence of a Riemannian
	metric of positive scalar curvature metric on manifolds with
	the Sullivan--Baas singularities. The manifolds we consider are
	$Spin$ and simply connected. We prove an analogue of the
	Gromov--Lawson Conjecture for such manifolds in the case of
	particular type of singularities. We give an affirmative
	answer when such manifolds with singularities accept a metric
	of positive scalar curvature in terms of the index of the
	Dirac operator valued in the corresponding ``$K$--theories with
	singularities''. The key ideas are based on the construction
	due to Stolz, some stable homotopy theory, and the index
	theory for the Dirac operator applied to the manifolds with
	singularities.  As a side-product we compute homotopy types of
	the corresponding classifying spectra.
\end{abstract}

\asciiabstract{We study the question of existence of a Riemannian
metric of positive scalar curvature metric on manifolds with the
Sullivan-Baas singularities. The manifolds we consider are Spin and
simply connected. We prove an analogue of the Gromov-Lawson Conjecture
for such manifolds in the case of particular type of singularities. We
give an affirmative answer when such manifolds with singularities
accept a metric of positive scalar curvature in terms of the index of
the Dirac operator valued in the corresponding "K-theories with
singularities". The key ideas are based on the construction due to
Stolz, some stable homotopy theory, and the index theory for the Dirac
operator applied to the manifolds with singularities.  As a
side-product we compute homotopy types of the corresponding
classifying spectra.}

\keywords{Positive scalar curvature, $Spin$ manifolds, manifolds
with singularities, $Spin$ cobordism, characteristic classes in
$K$--theory, cobordism with singularities, Dirac operator, $K$--theory
with singularities, Adams spectral sequence, $\A(1)$--modules.}

\asciikeywords{Positive scalar curvature, Spin manifolds, manifolds
with singularities, Spin cobordism, characteristic classes in
K-theory, cobordism with singularities, Dirac operator, K-theory
with singularities, Adams spectral sequence, A(1)-modules.}

\primaryclass{57R15}
\secondaryclass{53C21, 55T15, 57R90}
\maketitlepage

\begin{Relax}\end{Relax}
\let\arrow\undefined
\makeatletter
\input pb-diagram.sty
\makeatother

\section{Introduction}\label{s1}
{\bf \ref{s1}.1\qua Motivation}\qua It is well-known that the question of
existence of positive scalar curvature metric is hard enough for
regular manifolds. This question was studied extensively, and it is
completely understood, see \cite{GL1}, \cite{St1}, for simply
connected manifolds and for manifolds with few particular fundamental
groups, see \cite{BGS}, and also \cite{RS1}, \cite{RS2} for a detailed
discussion. At the same time, the central statement in this area, the
Gromov--Lawson--Rosenberg Conjecture is known to be false for some
particular manifolds, see \cite{Schick}. To motivate our interest we
first address a couple of naive questions. We shall consider manifolds
with boundary, and we always assume that a metric on a manifold is
product metric near its boundary. We use the abbreviation ``psc'' for
``positive scalar curvature'' throughout the paper.

Let $(P, g_P)$ be a closed Riemannian manifold, where the metric $g_P$
is not assumed to be of positive scalar curvature.  Let $X$ be a
closed manifold, such that the product $X\times P$ is a boundary of a
manifold $Y$. 

\medskip
{\bf Naive Question 1}\qua  {\sl Does there exist a psc-metric $g_X$ on
	$X$, so that the product metric $g_X\times g_P$ could be
	extended to a psc-metric $g_Y$ on $Y$?}

\medskip
\parbox{2.0in} {\noindent {\bf Examples}\qua (1)\qua Let $P=\<k\>=\{\mbox{$k$
	points}\}$, then a manifold $Y$ with $\p Y = X\times \<k\>$ is called
	a {\sl $\Z/k$--manifold.}  When $k=1$ (or $X=\p Y$) the above question
	is essentially trivial. Say, if $X$ and $Y$ are simply connected
	$Spin$ manifolds, and $\dim X=n-1\geq 5$, there is always a psc-metric
	$g_X$ which could be extended to a psc-metric $g_Y$.  }
\parbox{3.1in}{
\hspace*{4mm}\PSbox{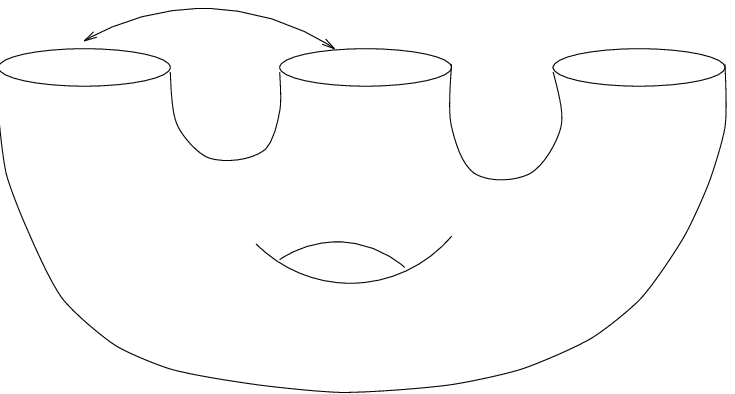}{5cm}{35mm}
\begin{picture}(10,0)
\put(-190,80){{\small $X$}}
\put(-110,80){{\small $X$}}
\put(-30,80){{\small $X$}}
\put(-110,20){{\small $Y$}}
\end{picture}
\centerline{{\small Figure 1: $\Z/k$--manifold}}}

	To see this one can delete a small open disk
	$D^n\subset Y$, and then push the standard metric on $S^{n-1}$
	through the cobordism $W=Y\setminus D^n$ to the manifold $X$
	using the surgery technique due to Gromov, Lawson \cite{GL1}
	and Schoen, Yau \cite{SY}.	

(2)\qua The case $P=\<k\>$ with $k\geq 2$ is not as simple. For
	example, there are many simply connected $Spin$ manifolds $X$
	of dimension $4k$ (for most $k$) which are not cobordant to
	zero, and, in the same time, two copies of $X$ are. Let $\p Y
	= 2 X$. It is not obvious that one can find a psc-metric $g_X$
	on $X$, so that the product metric $g_X\times \<2\>$ extends
	to a psc-metric $g_Y$ on $Y$.

(3)\qua Let $\Sigma^m$ (where $m=8l+1$ or $8l+2$, and $l\geq 1$) be
	a homotopy sphere which does not admit a psc-metric, see
	\cite{Hitchin}. We choose $k\geq 2 $ disjoint discs
	$D^{m}_1,\ldots, D^{m}_k\subset \Sigma^m$ and delete their
	interior. The resulting manifold $Y^m$ has the boundary
	$S^{m-1} \times \<k\>$. Clearly it is not possible to extend
	the standard metrics on the spheres $S^{m-1} \times \<k\>$ to
	a psc-metric on the manifold $Y$ since otherwise it would give
	a psc-metric on the original homotopy sphere
	$\Sigma^m$. However, it is not obvious that for {\it any
	choice of a psc-metric $g$ on $S^{m-1}$} the metric $g\times
	\<k\>$ could not be extended to a psc-metric on $Y^m$.

(4)\qua Let $P$ be again $k$ points. Consider a Joyce manifold $J^8$
	($Spin$, simply connected, Ricci flat, with $\hat{A}(J^8)=1$,
	and holonomy $Spin(7)$), see \cite{J}. Delete $k$ open disks
	$D^{m}_1,\ldots, D^{m}_k\subset J^8$ to obtain a manifold $M$,
	with $\p M = S^{7} \times \<k\>$. Let $g_0$ be the standard
	metric on $S^7$. Then clearly the metric $g_0\times \<k\>$ on
	the boundary $S^{7} \times \<k\>$ cannot be extended to a
	psc-metric on $M$ since otherwise one would construct a
	psc-metric on $J^8$. However, there are so called ``exotic''
	metrics on $S^7$ which are not in the same connective
	component as the standard metric. Nevertheless, as we shall
	see, there is no any psc-metric $g^{\prime}$ on $S^7$, so that
	the metric $g^{\prime}\times \<k\>$ could be extended to a
	psc-metric on $M$.

(5)\qua Let $P=S^1$ with nontrivial $Spin$ structure, so that $[P]$
	is a generator of the cobordism group $\Omega^{Spin}_1=\Z/2$.

\noindent
\parbox{2.8in}{\noindent
Let $d\theta^2$ be the standard metric
	on the circle. The analysis of the ring structure of
	$\Omega_*^{Spin}$ shows that there exist many examples of
	simply connected manifolds $X$ which are not $Spin$ cobordant
	to zero, however, the products $X\times P$ are, say $\p Y=
	X\times P$.

Again, in general situation there is no obvious clue whether for some
psc-metric $g_X$ on $X$ the product metric $g_X + d\theta^2$ on
$X\times P$ could be extended to a psc-metric on $Y$ or not.  }
\parbox{2.2in}{\noindent\vspace*{3mm}
\hspace*{6mm}
\PSboxa <-5pt, -5pt> {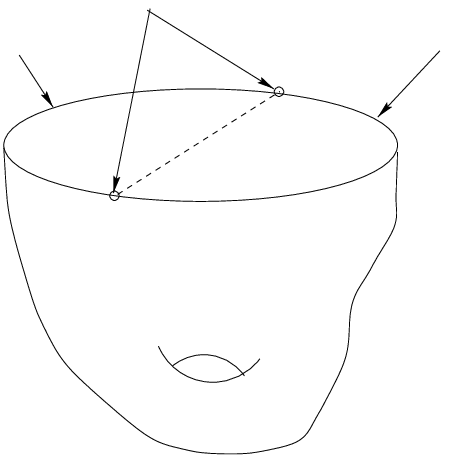}{5cm}{40mm}
\begin{picture}(0,0)
\put(30,150){{\small $X\times P_1\times P_2$}}
\put(10,135){{\small $Z_1\times P_1$}}
\put(120,137){{\small $Z_2\times P_2$}}
\put(70,60){{\small $Y$}}
\put(70,0){{\small Figure 2}}
\end{picture}
}

Now let $(P_1,g_1)$, $(P_2,g_2)$ be two closed Riemannian manifolds,
again, the metrics $g_1$, $g_2$ are not assumed to be of positive
scalar curvature. Let $X$ be a closed manifold such that
\begin{itemize}\leftskip -15pt
\item
the product $X\times P_1$ is a boundary of a manifold $Z_2$,
\item
the product $X\times P_2$ is a boundary of a manifold
$Z_1$,
\item the manifold $Z=Z_1\times P_1 \cup \epsilon
Z_2\times P_2$ is a boundary of a manifold $Y$ (where $\epsilon$ is an
appropriate sign if the manifolds are oriented), see Figure 2.
\end{itemize}

\medskip
{\bf Naive Question 2}\qua {\sl Does there exist a psc-metric $g_X$ on
	$X$, so that}

\noindent
(a)\qua {\sl the product metric $g_X\times g_1$ on $X\times P_1$ could
	be extended to a psc-metric $g_{Z_2}$ on $Z_2$,}

\noindent
(b)\qua {\sl the product metric $g_X\times g_2$ on $X\times P_1$ could
	be extended to a psc-metric on $g_{Z_1}$ $Z_1$,}

\noindent
(c)\qua {\sl the metric $g_{Z_1}\times g_{P_1}\cup g_{Z_2}\times
	g_{P_2}$ on the manifold $Z=Z_1\times P_1 \cup \epsilon
	Z_2\times P_2$ could be extended to a psc-metric $g_{Y}$ on
	$Y$?}

\medskip
{\bf \ref{s1}.2\qua Manifolds with singularities}\qua  Perhaps, one can
recognize that the above naive questions are actually about the
existence of a psc-metric on a manifold with the Baas--Sullivan
singularities, see \cite{S}, \cite{Baas}. In particular, a
$\Z/k$--manifold $M$ is a manifold with boundary $\p M$ diffeomorphic
to the product $\beta M\times \<k\>$. Then a metric $g$ on $M$ is a
regular Riemannian metric on $M$ such that it is product metric near
the boundary, and its restriction on each two components $\beta
M\times \{i\}$, $\beta M\times \{j\}$ are isometric via the above
diffeomorphism. To get the singularity one has to identify the
components $\beta M\times \{i\}$ with a single copy of $\beta M$.
Similarly a Riemannian metric may be defined for the case of general
singularities. We give details in Section \ref{s7}.

Thus manifolds with the Baas--Sullivan singularities provide an
adequate environment to reformulate the above naive question.  Let
$\Sigma =(P_1,\ldots,P_q)$ be a collection of closed manifolds, and
$M$ be a $\Sigma$--manifold (or manifold with singularities of the type
$\Sigma$), see \cite{Baas}, \cite{Mir}, \cite{Bot1} for definitions.
For example, if $\Sigma=(P)$, where $P=\<k\>$, a $\Sigma$--manifold $M$ is
$\Z/k$--manifold. Then the above questions lead to the following
one:

\medskip
{\bf Question}\qua {\sl Under which conditions does a $\Sigma$--manifold
	$M$ admit a psc-metric?}

\medskip
Probably it is hard to claim anything useful for a manifold with arbitrary
singularities. We restrict our attention to  $Spin$
simply connected manifolds and very particular singularities. 
Now we introduce necessary notation.

Let $\Omega^{Spin}_*(\cdot)$ be the $Spin$--cobordism theory, and
$MSpin$ be the Thom spectrum classifying this theory. Let
$\Omega^{Spin}_*(pt)=\Omega^{Spin}_*$ be the coefficient ring. Let
$P_1=\< 2 \> =\{\mbox{two points}\}$, $P_2$ be a circle with a
nontrivial $Spin$ structure, so that $[P_2]=\eta\in
\Omega^{Spin}_1\cong \Z/2$, and $P_3$, $[P_3]\in \Omega^{Spin}_8$, is
a Bott manifold, ie, a simply-connected manifold such that
$\widehat{A}(P_3) =1$.  There are different representatives of the
Bott manifold $P_3$.  Perhaps, the best choice is the Joyce manifold
$J^8$, \cite{J}.  Let $\Sigma_1=(P_1)$, $\Sigma_2=(P_1,P_2)$,
$\Sigma_3 = (P_1,P_2,P_3)$, and $\eta=(P_2)$.  We denote by
$\Omega^{Spin,\Sigma_i}_*(\cdot)$ the cobordism theory of
$Spin$--manifolds with $\Sigma_i$--singularities, and by
$MSpin^{\Sigma_i}$ the spectra classifying these theories, $i=1,2,3$.
We also study the theory $\Omega^{Spin,\eta}_*(\cdot)$, and the
classifying spectrum for this theory is denoted as $MSpin^{\eta}$. We
use notation $\Sigma$ for the above singularities $\Sigma_1$,
$\Sigma_2$, $\Sigma_3$ or $\eta$.

Let $KO_*(\cdot)$ be the periodic real $K$--theory, and $KO$ be the
classifying $\Omega$--spectrum. The Atiyah--Bott--Shapiro homomorphism
$\alpha \co \Omega^{Spin}_* \longrightarrow KO_*$ induces the map of spectra
\begin{equation}\label{eq1}
\alpha \co MSpin \longrightarrow  KO.
\end{equation}
It turns out that for our choice of singularities $\Sigma$ the
spectrum $MSpin^{\Sigma}$ splits as a smash product
$MSpin^{\Sigma}=MSpin\wedge X_{\Sigma}$ for some spectra $X_{\Sigma}$
(see Theorems \ref{Sigma1}, \ref{Sigma3}). We would like to introduce
the real $K$--theories $KO^{\Sigma}_*(\cdot)$ with the singularities
$\Sigma$. We define the classifying spectrum for
$KO^{\Sigma}_*(\cdot)$ by $ KO^{\Sigma} = KO\wedge X_{\Sigma}$. The
$K$--theories $KO^{\Sigma}_*(\cdot)$ may be identified with the
well-known $K$--theories. Indeed,
$$
KO^{\Sigma_1}_*(\cdot)= KO_*(\cdot;\Z/2), \ \ \ \
KO^{\eta}_*(\cdot)= K_*(\cdot), \ \ \ \ 
KO^{\Sigma_2}_*(\cdot)= K_*(\cdot;\Z/2),
$$
see Corollary \ref{KO-2}. The $K$--theory $KO_*^{\Sigma_3}(\cdot)$ is
``trivial'' since the classifying spectrum $KO^{\Sigma_3}$ is
contractible, see Corollary \ref{KO-3}.  Now the map $\alpha$ from
(\ref{eq1}) induces the map
$$
\alpha^{\Sigma} \co MSpin^{\Sigma}=MSpin\wedge X_{\Sigma}  
\stackrel{\alpha\wedge 1}{\longrightarrow}
KO\wedge X_{\Sigma} = KO^{\Sigma}
$$
and the homomorphism of the coefficient rings
\begin{equation}\label{eq2}
\alpha^{\Sigma} \co \Omega_*^{Spin,\Sigma} \longrightarrow KO_*^{\Sigma} .
\end{equation}
We define the integer $d(\Sigma)$ as follows:
$$
d(\Sigma_1)= 6, \ \ \ \ d(\Sigma_2) = 8, \ \ \ \ d(\Sigma_3)= 17,
\ \ \ \ d(\eta)= 7.
$$
Recall that if $M$ is a $\Sigma$--manifold, then (depending on the
length of $\Sigma$), the manifolds $\beta_i M$, $\beta_{ij}M$,
$\beta_{ijk} M$ (as $\Sigma$--manifolds) are defined in canonical
way. In particular, for $\Sigma=\Sigma_1,\eta$, there is a manifold
$\beta_i M$ such that $\p M = \beta_i M\times P_i$, for
$\Sigma=\Sigma_{2}$, there are $\Sigma$--manifolds $\beta_1 M$,
$\beta_2 M$,$\beta_{12} M$, and for $\Sigma=\Sigma_3$ there are
$\Sigma$--manifolds $\beta_i M$, $\beta_{ij}M$, $\beta_{ijk}M$.  These
manifolds may be empty.  The manifolds $\beta_i M$, $\beta_{ij} M $
and $\beta_{ijk} M$ are called $\Sigma$--{\sl strata of $M$.}

We say that a $\Sigma$--manifold $M$ is {\sl simply connected} if $M$
itself is simply connected and all $\Sigma$--strata of $M$ are simply
connected manifolds. 

\medskip
{\bf \ref{s1}.3\qua Main geometric result}\qua
The following theorem is the main geometric result of this paper.
\begin{Theorem}\label{ThA}
	Let $M^n$ be a simply connected $Spin$ $\Sigma$--manifold of
	dimension $n\geq d(\Sigma)$, so that all $\Sigma$--strata
	manifolds are nonempty manifolds. Then $M$ admits a metric of
	positive scalar curvature if and only if
	$\alpha^{\Sigma}([M])=0$ in the group $KO_n^{\Sigma}$.
\end{Theorem}
We complete the proof of Theorem \ref{ThA} only at the end of the
paper. However, we would like to present here the overview of the
main ingredients of the proof.

\medskip
{\bf \ref{s1}.4\qua Key ideas and constructions of the proof}\qua There are
two parts of Theorem \ref{ThA} to prove. The first ``if'' part is
almost ``pure topological''. The second ``only if'' part has more
analytical flavor.  We start with the topological ingredients.

The first key construction which allows to reduce the question on the
existence of a psc-metric to a topological problem, is the Surgery
Lemma. This fundamental observation originally is due to Gromov--Lawson
\cite{GL1} and Schoen--Yau \cite{SY}. We generalize the Surgery Lemma
for simply connected $Spin$ $\Sigma$--manifolds.

\noindent 
\parbox{2.5in} {\noindent 
This generalization is almost
straightforward, however we have to describe the surgery procedure for
$\Sigma$--manifolds.  

To explain the difference with the
case of regular surgery, we consider the example when $M$ is a
$\Z/k$--manifold, ie, $\p M = \beta M\times \<k\>$.  There are two
types of surgeries here. The first one is to do  surgery on
the interior of $M$, and the second one is to do surgery on
each manifold $\beta M$. 
}
\parbox{2.6in}{
\hspace*{2mm}\PSbox{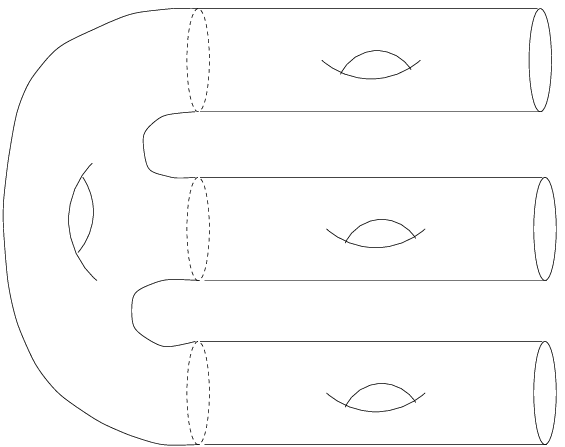}{5cm}{45mm}
\begin{picture}(0,0)
\put(-135,40){{\small $M$}}
\put(-70,110){{\small $V$}}
\put(-70,60){{\small $V$}}
\put(-70,10){{\small $V$}}
\put(15,110){{\small $\beta^{\prime}M$}}
\put(15,60){{\small $\beta^{\prime}M$}}
\put(15,10){{\small $\beta^{\prime}M$}}
\end{picture}
\vspace{3mm}

\centerline{{\small Figure 3: The manifold $M^{\prime}$}}}

We start with the second one.  Let $M$ be a $\Z/k$-manifold, with a
psc-metric $g_{M}$. We have $\p M = \beta M \times \<k\>$, where
$g_{\beta M}$ is a psc-metric.  Let $S^{p}\times D^{n-p-1}\subset\beta
M$, and $V$ be a trace of the surgery along the sphere $S^p$, ie, $\p
V = -\beta M \cup \beta^{\prime} M$.  We assume that $n-p-1\geq 3$, so
we can use the regular Surgery Lemma to push a psc-metric through the
manifold $V$ to obtain a psc-metric $g_{V}$ which is a product near
the boundary. Then we attach $k$ copies of $V$ to obtain a manifold $
M^{\prime}= M \cup_{\p M} V\times \<k\>, $ see Figure 3. Clearly the
metrics $g_M$ and $g_V$ match along a color of the common boundary,
giving a psc-metric $g^{\prime}$ on $M^{\prime}$.  

The first type of surgery is standard.  Let $S^{\ell}\times
D^{n-\ell}\subset M$ be a sphere together with a tubular neighborhood
inside the interior of the manifold $M$. Denote by $M^{\prime\prime}$
the result of surgery on $M$ along the sphere $S^{\ell}$. Notice that
$\p M^{\prime\prime}= \p M$. Then again the regular Surgery Lemma
delivers a psc-metric on $M^{\prime\prime}$.

The case of two and more singularities requires a bit more care. We
discuss the general Surgery procedure for $\Sigma$--manifolds in
Section \ref{s7}.  The Bordism Theorem (Theorem \ref{rev:Th3}) for
simply connected $\Sigma$--manifolds reduces the existence question of
a positive scalar curvature to finding a $\Sigma$--manifold within the
cobordism class $[M]_{\Sigma}$ equipped with a psc-metric.

To solve this problem we use the ideas and results due to S Stolz
\cite{St1}, \cite{St2}. The magic phenomenon discovered by S Stolz is
the following. Let us start with the quaternionic projective space
$\H\P^2$ equipped with the standard metric $g_0$ (of constant positive
curvature). It is not difficult to see that the  Lie group 
$$
G = PSp(3)=Sp(3)/\mbox{Center},
$$
acts by isometries of the metric $g_0$ on $\H\P^2$. Here
$\mbox{Center}\cong \Z/2$ is the center of the group $Sp(3)$. Then
given a smooth bundle $E \stackrel{p}{\longrightarrow}B$ of compact
$Spin$--manifolds, with a fiber $\H\P^2$, and a structure group $G$,
there is a straightforward construction of a psc-metric on the
manifold $E$, the total space of this bundle.  (A bundle with the
above properties is called a {\sl geometric $\H\P^2$--bundle}.) The
construction goes as follows. One picks an arbitrary metric $g_B$ on a
manifold $B$. Then locally, over an open set $U\subset B$, a metric on
$p^{-1}(U)\cong U\times \H\P^2$ is given as product metric
$g_E|_{p^{-1}(U)}=g_B|_U \times g_0$. By scaling the metric $g_0$, one
obtains that the scalar curvature of the metric $g_E|_{p^{-1}(U)}$ is
positive. Since the structure group of the bundle acts by isometries
of the metric $g_0$, one easily constructs a psc-metric $g_E$ on $E$.

Perhaps, this general construction was known for ages. The amazing
feature of geometric $\H\P^2$--bundles is that their total spaces, the
manifolds $E$, generate the kernel of the Atiyah--Bott--Shapiro
transformation $ \alpha \co \Omega^{Spin}_n \longrightarrow KO_n$. In
more detail, given an $\H\P^2$--bundle $E^n
\stackrel{p}{\longrightarrow} B^{n-8}$, there is a classifying map $f\co
B^{n-8} \longrightarrow BG$ which defines a cobordism class
$[(B,f)]\in \Omega^{Spin}_{n-8}(BG)$. The correspondence
$[(B,f)]\mapsto [E]\in \Omega^{Spin}_{n}$ defines the transfer map
$$
T \co \Omega^{Spin}_{n-8}(BG)  \longrightarrow \Omega^{Spin}_{n}.
$$
Stolz  proves \cite{St1} that $\mbox{Im}\ T = \mbox{Ker}\ \alpha$.   
Thus the manifolds $E$ deliver representatives in each cobordism class
of the kernel $\mbox{Ker}\ \alpha$.

We adopt this construction for manifolds with singularities. First we
notice that if a geometric $\H\P^2$--bundle $E \stackrel{p}{\longrightarrow}B$
is such that $B$ is a $\Sigma$--manifold, then $E$ is also a
$\Sigma$--manifold. In particular we obtain the induced transfer map
$$
T^{\Sigma} \co \Omega^{Spin,\Sigma}_*(BG) \longrightarrow
\Omega^{Spin,\Sigma}_{*+8}.
$$
The key here is to prove that $\mbox{Im}\ T^{\Sigma} = \mbox{Ker} \
\alpha^{\Sigma}$.  This requires complete information on the homotopy
type of the spectra $MSpin^{\Sigma}$. Sections \ref{s3}--\ref{s6} are
devoted to study of the spectra $MSpin^{\Sigma}$.

The second part, the proof of the ``only if'' statement, is geometric
and analytic by its nature.  We explain the main issues here for the
case of $\Z/k$--manifolds.  Recall that for a $Spin$ manifold $M$ the
direct image $\alpha([M]) \in KO_n$ is nothing else but the {\it
topological index} of $M$ which coincides (via the Atiyah--Singer index
theorem) with the {\it analytical index} $\ind(M)\in KO_n$ of the
corresponding Dirac operator on $M$. Then the Lichnerowicz formula and
its modern versions imply that the analytical index $\ind(M)$ vanishes
if there is a psc-metric on $M$.

Thus if we would like to give a similar line of arguments for
$\Z/k$--manifolds, we face the following issues.  To begin with, we
should have the Dirac operator to be well-defined on a $Spin$
$\Z/k$--manifold. Then we have to define the $\Z/k$--version of the
analytical index $\ind_{\Z/k}(M)\in KO_n^{\<k\>}$ and to prove the
vanishing result, ie, that $\ind_{\Z/k}(M)=0$ provided that there is
a psc-metric on $M$.  Thirdly we must identify the analytical index
$\ind_{\Z/k}(M)$ with the direct image $\alpha^{\<k\>}([M]) \in
KO_n^{\<k\>}$, ie, to prove the $\Z/k$--mod version of the index
theorem.  These issues were already addressed, and, in the case of
$Spin^c$--manifolds, resolved by Freed \cite{Freed}, \cite{Fr2}, Freed
\& Melrose \cite{FrMel}, Higson \cite{Higson}, Kaminker \&
Wojciechowski \cite{Kam-Woj}, and Zhang \cite{Zh1,Zh2}. Unfortunately,
the above papers study mostly the case of $Spin^c$ $\Z/k$--manifolds
(with the exception of \cite{Zh1,Zh2} where the mod 2 index is
considered), and the general case of $Spin$ $\Z/k$--manifolds is
essentially left out in the cited work.  The paper \cite{R} by
J. Rosenberg shows that the Dirac operator and its index are
well-defined for $\Z/k$--manifolds and there the index vanishes if a
$Spin$ $\Z/k$--manifold has psc-metric.  The case of general
singularities $\Sigma$ require more work. Here we use the results of
\cite{R} to prove that if a $\Sigma$--manifold $M$ has a psc-metric,
then $\alpha^{\Sigma}([M])= 0$ in the group $KO^{\Sigma}$.  In order
to prove this fact we essentially use the specific homotopy features
of the spectra $MSpin^{\Sigma}$.

The plan is the following. We give necessary definitions and
constructions on manifolds with singularities in Section \ref{s2}.
The next four sections are devoted to homotopy-theoretical study of
the spectra $MSpin^{\Sigma}$. We describe the homotopy type of the
spectra $MSpin^{\Sigma_1}$, $MSpin^{\Sigma_2}$, and $MSpin^{\eta}$ in
Section \ref{s3}. We describe a product structure of these spectra in
Section \ref{s4}. In Section \ref{s5} we describe a splitting of the
spectra $MSpin^{\Sigma}$ into indecomposable spectra. In Section
\ref{s6} we describe the homotopy type of the spectrum
$MSpin^{\Sigma_3}$.  We prove the Surgery Lemma for manifolds with
singularities in Section \ref{s7}.  Section \ref{s8} is devoted to the
proof of Theorem \ref{ThA}.

It is a pleasure to thank Hal Sadofsky for helpful discussions on the
homotopy theory involved in this paper, and acknowledge my
appreciation to Stephan Stolz for numerous discussions about the
positive scalar curvature.  The author also would like to thank the
Department of Mathematics of the National University of Singapore for
hospitality (this was Fall of 1999). The author is thankful to
Jonathan Rosenberg for his interest to this work and useful
discussions. Finally, the author thanks the referee for helpful
suggestions.

\section{Manifolds with singularities}\label{s2}
Here we briefly recall basic definitions concerning manifolds with the
Baas--Sullivan singularities. Let $G$ be a stable Lie group. We will be
interested in the case when $G=Spin$. Consider the category of smooth
compact manifolds with a stable $G$--structure in their stable normal
bundle.

\medskip
{\bf \ref{s2}.1\qua General definition}\qua 
Let $\Sigma= (P_1,\ldots,P_k)$, where $P_1,\ldots,P_k$ are arbitrary
closed manifolds (possibly empty). It is convenient to denote
$P_0=pt$. Let $I = \{i_1,\ldots ,i_q\} \subset \{0,1,\ldots ,k\}$. We
denote $ P^{I} = P_{i_1}\X \ldots \X P_{i_q}$.

\begin{Definition} 
	{\rm We call a manifold $M$ a $\Sigma$--{\sl manifold} if there
	are given the following:

(i)\qua a partition $ \p M = \p_0M \cup \p_1M \cup \ldots \cup \p_kM
	$ of its boundary $\p M$ such that the intersection $ \p_{I}M
	= \p_{i_1}M \cap \ldots \cap \p_{i_q}M $ is a manifold for
	every collection $I = \{i_1,\ldots ,i_q\} \subset \{0,1,\ldots
	,k\}$, and its boundary is equal to
$$
\p \( \p_{I}M \) =\bigcup_{j\notin I} \( \p_{I} M \cap
\p_j M \);
$$ 
(ii)\qua compatible product structures (ie, diffeomorphisms
	preserving the stable $G$--structure)
$$
\phi_{I} \co \p_{I} M \longrightarrow  \beta_{I} M \X  P^{I} .
$$
	Compatibility means that if $I \subset J$ and $ \iota \co \p_{J}
	M \longrightarrow \p_{I} M $ is the inclusion, then the map
$$
\phi_{I} \circ \iota \circ \phi^{-1}_{J} \co \beta_{J}M \X P^{J}
\longrightarrow \beta_{I} M \X P^{I} 
$$
	is identical on the direct factor $P^{I}$.}
\end{Definition}
To get actual singularities we do the following.  Two points $x, y$ of
a $\S$--manifold $M$ are {\sl equivalent} if they belong to the same
manifold $\p_{I}M$ for some $I \subset \{0, 1,\ldots ,k \}$ and $
pr\circ \phi_{I}(x) = pr\circ \phi_{J}(y), $ where $ pr \co \beta_{I}M
\X P^{I} \longrightarrow \beta_{I}M $ is the projection on the direct
factor. The factor-space of $M$ under this equivalence relation is
called {\sl the model of the} $\S$--{\sl manifold} $M$ and is denoted
by $M_{\S}$. Actually it is convenient to deal with $\S$--manifolds
without considering their models. Indeed, we only have to make sure
that all constructions are consistent with the projections $ \pi \co M
\longrightarrow M_{\S } .  $ The {\sl boundary} $\delta M$ {\sl of a}
$\S$--{\sl manifold} $M$ is the manifold $\p_0 M$. If $\delta
M=\emptyset$, we call $M$ a {\sl closed $\S$--manifold}.  The boundary
$\delta M$ is also a $\S$--manifold with the inherited decomposition $
\p_I (\delta M) = \p_{I} M \cap \delta M$. The manifolds $\beta_{I}M$
also inherit a structure of a $\S$--manifold:
\begin{equation}
\p_j (\beta_{I} M ) = \{\begin{array}{cl} 
\strut0{10}\emptyset & \mbox{if $j \in I$,} \\
\beta_{\{j\}\cup I} M \X P_j & \mbox{otherwise.}
		\end{array} \right.  
\end{equation}
Here we denote $ \beta_{I} M = \beta_{i_1} \( \beta_{i_2} \( \cdots
\beta_{i_q} M \) \cdots \)$ for $I = \{ i_1,\ldots ,i_q \} \subset \{
1,\ldots ,k \}$.

Let $(X,Y)$ be a pair of spaces, and $f\co (M,\delta M) \longrightarrow
(X,Y)$ be a map. Then the pair $(M,f)$ is a {\sl singular
$\Sigma$--manifold} of $(X,Y)$ if the map $f$ is such that for every
index subset $I = \{ i_1,\ldots ,i_q \} \subset \{ 1,\ldots ,k \}$ the
map $f|_{\p_{I}M} $ is decomposed as $ f|_{\p_{I}M} = f_{I}\circ
pr\circ \phi_I $, where the map $\phi_I$ as above, $pr \co \beta_I M
\times P^I \longrightarrow \beta_I M $ is the projection on the direct
factor, and $f_I \co \beta_I M \longrightarrow X$ is a continuous
map. The maps $f_I$ should be compatible for different indices $I$ in
the obvious sense.
\begin{Remark} Let $(M,f)$ be a {\sl singular $\Sigma$--manifold}, then
	the map $f$ factors through as $f= f_{\Sigma}\circ \pi$, where
	$\pi \co M \longrightarrow M_{\Sigma}$ is the canonical
	projection, and $f_{\Sigma} \co M_{\Sigma} \longrightarrow X$ is
	a continuous map.  We also notice that singular
	$\Sigma$--manifolds may be identified with their
	$\Sigma$--models.
\end{Remark}
The cobordism theory $\Omega_*^{G,\Sigma}(\cdot)$ of
$\Sigma$--manifolds is defined in the standard way. In the case of
interest, when $G=Spin$, we denote $MSpin^{\Sigma}$ a spectrum
classifying the cobordism theory $\Omega_*^{Spin,\Sigma}(\cdot)$.

\medskip
{\bf \ref{s2}.2\qua The case of two and three singularities}\qua We start
with the case $\Sigma=(P_1,P_2)$. 
Then if $M$ is a $\Sigma$--manifold,
we have that the diffeomorphisms
$$
\begin{array}{c}
\strut0{10}\phi \co \p M \stackrel{\cong}{\longrightarrow} \p_1 M \cup \p_2 M, 
\\ 
\strut0{10}\phi_i \co \p_i M \stackrel{\cong}{\longrightarrow} \beta_i M \times
P_i, \ \ \ i=1,2; 
\\ 
\phi_{12} \co \p_1 M \cap \p_2 M
\stackrel{\cong}{\longrightarrow} \beta_{12} M \times P_1\times P_2
\end{array}
$$
are given. We always assume that the manifold 
$\beta_{12} M \times
P_1\times P_2$ is embedded into  $\p_1 M$ and $\p_2 M $ together with a color:

\noindent
\parbox{3.0in} {\noindent 
$$
\beta_{12} M \times P_1\times P_2 \times I \subset \p_1 M, \p_2 M.
$$
Thus we actually have the following decomposition of the boundary $\p M$:
$$
\p M \cong \p_1 M \cup \(\beta_{12}M \times P_1 \times P_2\times I\)\cup
\p_2 M,
$$
so the manifold $\beta_{12}M \times P_1 \times P_2$ is ``fattened'' inside
$\p M$. Also we assume that the boundary $\p M$ is embedded into $M$
together with a color $\p M \times I\subset M$, see Figure 4.
}\hfill
\parbox{2in}{
\hspace*{5mm}\PSboxa <0pt,-10pt> {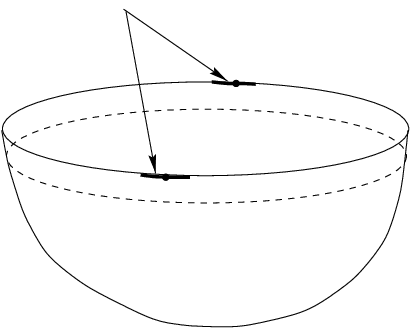}{5cm}{38mm}
\begin{picture}(0,0)
\put(65,20){{\small $M$}}
\put(10,100){{\small $\beta_{12}M\!\times \!P_1\!\times \!P_2\times \! I$}}
\put(120,73){{\small $\p_1 M$}}
\put(-7,63){{\small $\p_2 M$}}
\end{picture}

\centerline{{\small Figure 4}}}

The case when  $\Sigma=(P_1,P_2,P_3)$ is the
most complicated one we are going to work with. 

\noindent
\parbox{2.6in} {\noindent 
Let $M$ be a
closed $\Sigma$--manifold, then we are given the diffeomorphisms:
$$
\begin{array}{l}
\strut0{10}\phi \co \p M \stackrel{\cong}{\longrightarrow} \p_1 M \cup \p_2 M \cup \p_3 M,
\\
\strut0{10}\phi_i \co \p_i M \stackrel{\cong}{\longrightarrow} \beta_i M 
\times P_i, \ \ \ i=1,2,3;
\\
\strut0{10}\phi_{ij} \co  \p_i M \cap \p_j M \stackrel{\cong}{\longrightarrow} 
\beta_{ij} M \times P_i\times P_j,
\\
\strut0{10}\phi_{123} \co  \p_1 M \cap \p_2 M \cap \p_3 M 
\stackrel{\cong}{\longrightarrow} 
\\
\mbox{ \ } \ \ \ \beta_{123} M \times P_1\times P_2\times P_3
\end{array}
$$
where $i,j =1,2,3, \ \ i\neq j$, see Figure 5.  
}
\parbox{2.2in}{
\hspace*{3mm}\PSbox{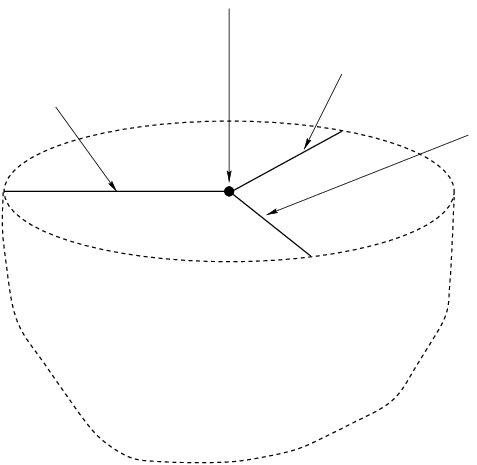}{5cm}{55mm}
\begin{picture}(0,0)
\put(-85,40){{\small $M$}}
\put(-75,118){{\small $\beta_{12}M\!\times \!P_1\!\times \!P_2$}}
\put(-40,100){{\small $\beta_{13}M\!\times \!P_1\!\times \!P_3$}}
\put(-160,105){{\small $\beta_{23}M\!\times \!P_2\!\times \!P_3$}}
\put(-135,140){{\small $\beta_{123}M\!\times \!P_1\!\times \!P_2\!\times \!P_3$}}
\put(-50,70){{\small $\p_1 M$}}
\put(-110,85){{\small $\p_2 M$}}
\put(-110,65){{\small $\p_3 M$}}
\end{picture}

\centerline{{\small Figure 5}}}

First, we assume here that the boundary $\p M$ is embedded into $M$
together with a color $(0,1]\times \p M$.  The decomposition 
$$
\p M
\stackrel{\phi}{\longrightarrow} \p_1 M \cup \p_2 M \cup \p_3 M
$$
gives also the ``color''
structure on $\p M$. 

We assume that the boundary $\p(\p_i M)$ is
embedded into $\p_i M$ together with the color 
$(0,1]\times \p(\p_i
M).
$

Even more, we assume that the manifold 
$
\beta_{123} M \times P_1\times P_2\times P_3
$ 
is embedded into the boundary $\p M$ together
with its normal tube:
$$
\beta_{123} M \times P_1\times P_2\times P_3 \times D^2 \subset \p M,
$$

\parbox{57mm}{so that the colors of the manifolds 
$$
\beta_{ij} M \times P_i \times P_j\subset \p_iM\cap \p_j M
$$
are compatible with this embedding, as is shown on Figure 6. 
As in the case of two singularities, the submanifolds 
$$
\begin{array}{l}
\strut0{10}\beta_{ij}M \times P_i \times P_j \ \ \ 
\mbox{and} \ \ \
\\
 \beta_{123}M \times P_1 \times P_2 \times P_3
\end{array}
$$
are ``fattened'' inside the boundary $\p M$.  Furthermore, we assume
that there are not any corners in the above color decomposition.
}
\noindent
\parbox{2.9in}{
\hspace*{8mm}
\PSbox{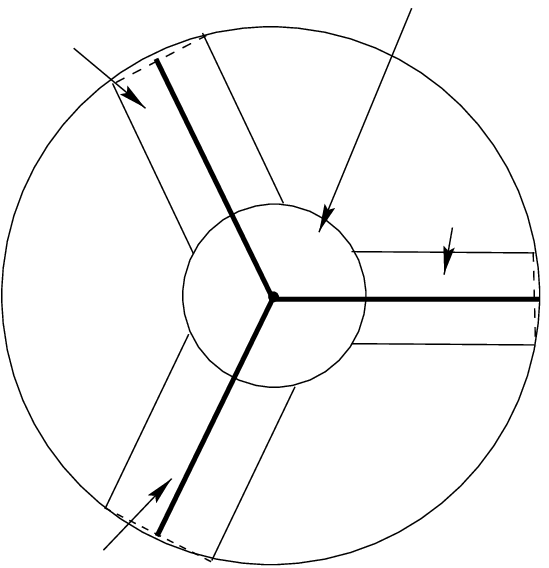}{5cm}{60mm}
\begin{picture}(0,0)
\put(-45,100){{\small $\beta_{12}M\!\times \!P_1\!\times \!P_2\!\times \!I$}}
\put(-160,-10){{\small $\beta_{13}M\!\times \!P_1\!\times \!P_3\!\times \!I$}}
\put(-170,158){{\small $\beta_{23}M\!\times \!P_2\!\times \!P_3\!\times \!I$}}
\put(-75,165){{\small $\beta_{123}M\!\times \!P_1\!\times
\!P_2\!\times \!P_3\!\times \! D^2$}} 
\put(-50,30){{\small $\p_1 M$}}
\put(-70,125){{\small $\p_2 M$}} 
\put(-120,75){{\small $\p_3 M$}}
\end{picture}
\vspace*{4mm}

\centerline{{\small Figure 6}}}

\medskip
{\bf \ref{s2}.3\qua Bockstein--Sullivan exact sequence}\qua  Let $MG$ be the
Thom spectrum classifying the cobordism theory $\Omega^{G}_*(\cdots)$.
Let $\Sigma=(P)$, and $p=\dim P$. Then there is a stable map 
$S^p \stackrel{[P]}{\longrightarrow}
MG$ representing the element $[P]$. Then we have the composition
$$
\cdot [P]\co \Sigma^p MG=S^p \wedge MG
\stackrel{[P]\wedge Id}{\longrightarrow} MG \wedge MG
\stackrel{\mu}{\longrightarrow} MG
$$
where $\mu$ is the map giving $MG$ a structure of a ring spectrum.
Then the cofiber, the spectrum $MG^{\Sigma}$ of the map
\begin{equation}\label{Bok1}
\Sigma^p MG \stackrel{\cdot [P]}{\longrightarrow} MG 
\stackrel{\pi}{\longrightarrow} MG^{\Sigma}
\end{equation}
is a classifying spectrum for the cobordism theory $\Omega_*^{G,\Sigma}$.
The cofiber (\ref{Bok1}) induce the long exact Bockstein--Sullivan sequence 
\begin{equation}\label{Bok2}
\!\!\!\!\!\cdots \rightarrow  \Omega_{n-p}^{G}(X,A)
\stackrel{\cdot [P]}{\longrightarrow} 
\Omega_{}^{G}(X,A) \stackrel{\pi}{\longrightarrow} 
\Omega_{n-p}^{G,\Sigma}(X,A) \stackrel{\beta}{\longrightarrow}
\Omega_{n-p-1}^{G}(X,A)\rightarrow \cdots
\end{equation}
for any $CW$--pair $(X,A)$.
Similarly, if $\Sigma_{j}=(P_1,\ldots,P_j)$, $j=1,\ldots,k$, then there
is a cofiber
$$
\Sigma^{p_j} MG^{\Sigma_{j-1}} \stackrel{\cdot [P_j]}{\longrightarrow}
MG^{\Sigma_{j-1}}
\stackrel{\pi_j}{\longrightarrow} MG^{\Sigma_j}
$$
induce the exact Bockstein--Sullivan sequence
\begin{equation}\label{Bok3}
\!\!\!\!\!\cdots 
\stackrel{\beta_j}{\longrightarrow}
\Omega_{n-p_j}^{G,\Sigma_{j-1}}(X,A) 
\stackrel{\cdot [P_j]}{\longrightarrow}
\Omega_{n}^{G,\Sigma_{j-1}}(X,A) 
\stackrel{\pi_j}{\longrightarrow}  
\Omega_{n}^{G,\Sigma_j}(X,A) 
\stackrel{\beta_j}{\longrightarrow}
\cdots\!\!\!\!\!
\end{equation}
for any $CW$--pair $(X,A)$. We shall use the Bockstein--Sullivan exact sequences 
(\ref{Bok2}), (\ref{Bok3})
throughout the paper.

\section{The spectra $MSpin^{\Sigma_1}$, $MSpin^{\Sigma_2}$ and
 $MSpin^{\eta}$}\label{s3} Let $M(2)$ be the mod 2 Moore spectrum with
 the bottom cell in zero dimension, ie, $M(2) =\Sigma^{-1}\R\P^2$. We
 consider also the spectrum $\Sigma^{-2}\C\P^2$ and the spectrum
 $Y=M(2)\wedge \Sigma^{-2}\C\P^2$ which was first studied by M Mahowald,
 \cite{Mah1}.  Here is the result on the spectra $MSpin^{\Sigma_1}$,
 $MSpin^{\Sigma_2}$ and $MSpin^{\eta}$.
\begin{Theorem}\label{Sigma1} There are homotopy
equivalences:
\begin{enumerate}
\item[{\rm (i)}] 
$MSpin^{\Sigma_1} \cong MSpin\wedge M(2)$,
\item[{\rm (ii)}] 
$MSpin^{\eta} \ \ \cong MSpin\wedge \Sigma^{-2}\C\P^2$,
\item[{\rm (iii)}]
$MSpin^{\Sigma_2} \cong MSpin\wedge Y$.
\end{enumerate}
\end{Theorem}
\begin{proof} Let $\iota \co S^0 \longrightarrow MSpin$ be a unit map.
The main reason why the above homotopy equivalences hold is that the
elements $2,\eta\in \Omega^{Spin}_*$ are in the image of the
homomorphism $\iota_* \co S_*^0 \longrightarrow
\Omega^{Spin}_*$. Indeed, consider first the spectrum
$MSpin^{\eta}$. Let $S^1 \stackrel{\eta}{\longrightarrow} S^0$ be a
map representing $\eta\in \pi_1(S^0)$.  We obtain the cofibration:
\begin{equation}\label{eta}
S^1 \stackrel{\eta}{\longrightarrow} S^0 
\stackrel{\pi}{\longrightarrow} \Sigma^{-2}\C\P^2 \ .
\end{equation}
Then the composition $S^1  \stackrel{\eta}{\longrightarrow} S^0 
 \stackrel{\iota}{\longrightarrow} MSpin$ represents
$\eta\in MSpin_1$. Let $\cdot\eta$ be the map
$$
\cdot\eta \co S^1\wedge MSpin \stackrel{\iota\eta\wedge 1}{\longrightarrow} 
MSpin\wedge MSpin \stackrel{\mu}{\longrightarrow}  MSpin,
$$
where $\mu$ is a multiplication. Note that the diagram
$$
\begin{CD}
S^1\wedge MSpin @>\iota\eta\wedge 1>> MSpin\wedge MSpin @>\mu>> MSpin 
\\
@A1\wedge 1 AA @A\iota\wedge 1 AA @A 1 AA 
\\
S^1 @>\eta\wedge 1>> S^0\wedge MSpin @>\cong>> MSpin 
\end{CD}
$$
commutes since the map $\iota \co S^0 \longrightarrow MSpin$ represents a unit of the ring spectrum $MSpin$. We obtain
a commutative diagram of cofibrations:
\begin{equation}\label{eta2}
\begin{CD}
S^1 \wedge MSpin @>\cdot\eta>> MSpin @>\pi_{\eta}>> MSpin^{\eta}
\\
@A1\wedge 1 AA @A 1 AA @A f_{\eta} AA
\\
S^1 \wedge MSpin @>\eta\wedge 1>> MSpin @>\pi\wedge 1>>
\Sigma^{-2}\C\P^2 \wedge MSpin
\end{CD}
\end{equation}
where $f_{\eta} \co MSpin^{\eta}\longrightarrow \Sigma^{-2}\C\P^2 \wedge
MSpin \cong MSpin \wedge\Sigma^{-2}\C\P^2$ gives a homotopy
equivalence by 5--lemma.  The proof for the spectrum
$MSpin^{\Sigma_1}=MSpin^{\<2\>} $ is similar.

Consider the spectrum $MSpin^{\Sigma_2}$. First we note that 
the bordism theory
$\Omega^{Spin,\Sigma_2}_*(\cdot)=\Omega^{Spin,(P_1,P_2)}_*(\cdot)$
coincides with the theory $\Omega^{Spin,(P_2,P_1)}_*(\cdot)$, where the
order of singularities is switched. In particular, the spectrum
$MSpin^{\Sigma_2}$ is a cofiber in the following
cofibration:
\begin{equation}\label{eta3}
\begin{CD}
S^0 \wedge MSpin^{\eta} @>\cdot 2 >> MSpin^{\eta} @>>>
MSpin^{\Sigma_2}.
\end{CD}
\end{equation}
Here the map $\cdot 2 \co S^0 \wedge MSpin^{\eta} \longrightarrow
MSpin^{\eta}$ is defined as follows. Let $S^0
\stackrel{2}{\longrightarrow} S^0$ be a map of degree 2.  Then the
composition $S^0 \stackrel{2}{\longrightarrow} S^0
\stackrel{\iota}{\longrightarrow} MSpin$ represents $2\in
\Omega^{Spin}_0$. The spectrum $MSpin^{\eta}$ is a module (say, left)
spectrum over $MSpin$, ie, there is a map $\mu^{\prime}_L :
MSpin\wedge MSpin^{\eta} \longrightarrow MSpin^{\eta}$ so that the diagram
$$
\begin{CD}
MSpin\wedge MSpin @>\mu>> MSpin 
\\
@V 1\wedge\pi_{\eta}VV @V\pi_{\eta}VV
\\
MSpin\wedge MSpin^{\eta} @>\mu^{\prime}_L>> MSpin^{\eta}
\end{CD}
$$
commutes. Then the map $\cdot 2$ is defined as composition:
$$
\begin{CD}
S^0\wedge MSpin^{\eta} @>2\iota \wedge 1>> MSpin\wedge MSpin^{\eta} 
@>\mu^{\prime}_L>> MSpin^{\eta}.
\end{CD}
$$
Note that the diagram
$$
\begin{CD}
S^0\wedge MSpin^{\eta} @>2\iota \wedge 1>> MSpin\wedge MSpin^{\eta} 
@>\mu^{\prime}_L>> MSpin^{\eta} 
\\
@A1\wedge 1 AA @A\iota\wedge 1 AA @A 1 AA 
\\
S^0\wedge  MSpin^{\eta} @> 2\wedge 1>> S^0\wedge MSpin^{\eta} @>\cong>> MSpin^{\eta} 
\end{CD}
$$
commutes since $S^0 \stackrel{\iota}{\longrightarrow} MSpin$
represents a unit, and $MSpin^{\eta}$ is a left module over the ring
spectrum $MSpin$. We obtain the commutative diagram of cofibrations:
\begin{equation}\label{eta4}
\begin{CD}
S^0 \wedge MSpin^{\eta} @>\cdot 2>> MSpin @>\pi_{2}>> MSpin^{\Sigma_2}
\\
@A1\wedge 1 AA @A 1 AA @A f_{2} AA
\\
S^0 \wedge MSpin^{\eta} @>2\wedge 1>> MSpin^{\eta} @>\pi\wedge 1>>
M(2)\wedge MSpin^{\eta}
\end{CD}
\end{equation}
The map $f_{2} \co M(2)\wedge MSpin^{\eta} \longrightarrow
MSpin^{\Sigma_2}$ gives a desired homotopy equivalence. Thus we have
$MSpin^{\Sigma_2}\cong M(2)\wedge MSpin^{\eta} \cong
MSpin^{\eta}\wedge M(2) =MSpin \wedge Y$.\hfill 
\end{proof}
\begin{Remark} In the above proof, we did not use any
specific properties of the spectrum $MSpin$ except that it is a ring
spectrum. In fact, $MSpin$ may be replaced by any other classic Thom
spectrum.
\end{Remark}
Later we prove that the homotopy equivalence
$$
MSpin^{\Sigma_3}
\sim MSpin \wedge \Sigma^{-2} \C\P^2\wedge V(1) \ ,
$$
where $V(1)$ is the cofiber of the Adams map $A \co \Sigma^{8}M(2)
\longrightarrow M(2)$. However, first we have to study the spectra
$MSpin^{\Sigma_1}$, $MSpin^{\Sigma_2}$ and $MSpin^{\eta}$ in more
detail.
\section{Product structure}\label{s4}
Recall that the spectrum $MSpin$ is a ring spectrum. Here we work with
the category of spectra, and commutativity of diagrams mean
commutativity up to homotopy. Let, as above, $\iota \co S^0
\longrightarrow MSpin$ be the unit, and $\mu \co MSpin \wedge MSpin
\longrightarrow MSpin$ the map defining the product structure. Let
$MSpin^{\Sigma}$ be one of the spectrum we considered above. The
natural map $\pi \co MSpin \longrightarrow MSpin^{\Sigma}$ turns the
spectrum $MSpin^{\Sigma}$ into a left and a right module over the
spectrum $MSpin$, ie, there are maps
$$
\begin{array}{c}
\mu^{\prime}_L\co MSpin \wedge  MSpin^{\Sigma} 
\longrightarrow  MSpin^{\Sigma},
\ \
\mu^{\prime}_R\co MSpin^{\Sigma} \wedge MSpin 
\longrightarrow  MSpin^{\Sigma},
\end{array}
$$
so that the diagrams
$$
\begin{CD} 
MSpin \wedge MSpin @>\mu>> MSpin
\\
@V 1\wedge \pi VV @V\pi VV
\\
MSpin \wedge MSpin^{\Sigma}@>\mu^{\prime}_L>> MSpin^{\Sigma}
\end{CD}\ \ \ \ 
\begin{CD}
MSpin \wedge MSpin @>\mu>> MSpin
\\
@V\pi \wedge 1 VV @V\pi VV
\\
MSpin^{\Sigma}  \wedge MSpin @>\mu^{\prime}_R>> MSpin^{\Sigma}
\end{CD}
$$
commute. We say that the spectrum
$MSpin^{\Sigma}$ has an {\sl admissible ring structure} 
$$
\mu^{\Sigma}\co MSpin^{\Sigma} \wedge MSpin^{\Sigma} \longrightarrow  MSpin^{\Sigma}
$$
if the map $S^0 \stackrel{\iota}{\longrightarrow} 
MSpin \stackrel{\pi}{\longrightarrow} MSpin^{\Sigma}$ is a unit,
and the diagrams
$$
\begin{CD}
MSpin \wedge MSpin^{\Sigma}@>\mu^{\prime}_L>> MSpin^{\Sigma}
\\
@V \pi\wedge 1 VV @V 1 VV
\\
MSpin^{\Sigma} \wedge MSpin^{\Sigma} @>\mu^{\Sigma}>> MSpin^{\Sigma}
\end{CD} \ \ \ \
\begin{CD}
MSpin^{\Sigma} \wedge MSpin  @>\mu^{\prime}_R>> MSpin^{\Sigma}
\\
@V 1 \wedge \pi VV @V 1 VV
\\
MSpin^{\Sigma} \wedge MSpin^{\Sigma} @>\mu^{\Sigma}>> MSpin^{\Sigma}
\end{CD} 
$$
commute. The questions of existence, commutativity and associativity
of an admissible product structure were thoroughly studied in
\cite{Bot1}, \cite{Mir}.
\begin{Theorem}\label{product} 
\begin{enumerate}
\item[{\rm (i)}] The spectrum $MSpin^{\Sigma_1}$
does not admit an admissible product structure.
\item[{\rm (ii)}] The spectra $MSpin^{\eta}$, $MSpin^{\Sigma_2}$ and
$MSpin^{\Sigma_3}$ have admissible product structures $\mu^{\eta}$,
$\mu^{\Sigma_2}=\mu^{(2)}$, and $\mu^{\Sigma_3}=\mu^{(3)}$  respectively.
\item[{\rm (iii)}] For any choice of an admissible product structure
$\mu^{\eta}$, it is commutative and associative.  For any choice of
admissible product structures $\mu^{(2)}$, and $\mu^{(3)}$,
they are associative, but not commutative.
\end{enumerate}
\end{Theorem}
\begin{proof}
Recall that for each singularity manifold $P_i$ there is
an obstruction manifold $P_i^{\prime}$ with singularity.  In the cases
of interest, we have: $[P_1^{\prime}]_{\Sigma_1} = \eta \in
\Omega^{Spin,\Sigma_1}_1$, which is non-trivial; and the obstruction
$[P_2^{\prime}]\in \Omega^{Spin,\Sigma_2}_3 =0$, and
$[P_2^{\prime}]\in \Omega^{Spin,\eta}_3 =0$. Thus \cite[Lemma
2.2.1]{Bot1} implies that there is no admissible product structure in
the cobordism theory $\Omega^{Spin,\Sigma_1}_*(\cdot)$, so the
spectrum $MSpin^{\Sigma_1}$ does not admit an admissible product
structure.  The obstruction element $[P_3^{\prime}]_{\Sigma_3}\in
\Omega^{Spin,\Sigma_3}_{17}$, and since $\dim P_3 =8$ is even, the
obstruction manifold $P_3^{\prime}$ is, in fact, a manifold without
any singularities (see \cite{Mir}), so the element
$[P_3^{\prime}]_{\Sigma_3}$ is in the image $\Im (\Omega^{Spin}_{17}
\longrightarrow  \Omega^{Spin,\Sigma_3}_{17})$.  However, the elements of
$\Omega^{Spin}_{17}$ are divisible by $\eta$, so they are zero in the
group $\Omega^{Spin,\eta}_{17}$, and, consequently, in
$\Omega^{Spin,\Sigma_3}_{17}$.

The result of \cite[Theorem 2.2.2]{Bot1} implies that the spectra
$MSpin^{\eta}$, $MSpin^{\Sigma_2}$ and $MSpin^{\Sigma_3}$ have
admissible product structures $\mu^{(2)}$ and $\mu^{\eta}$
respectively.

It is also well-known \cite{Wurg} that the element $v_1\in
\Omega^{Spin,\Sigma_2}_2$ is an obstruction to the commutativity of
the product structure $\mu^{(2)}$. An obstruction to the commutativity
for the product structure $\mu^{\eta}$ lives in the group
$\Omega^{Spin,\Sigma_2}_5=0$.  The obstructions to associativity are
3--torsion elements, (see \cite[Lemma 4.2.4]{Bot1}) so they all are
zero. \hfill 
\end{proof}
\section{Homotopy structure of the spectra $MSpin^{\Sigma}$}\label{s5}
First we recall the work of Anderson, Brown, and Peterson \cite{ABP}
on structure of the spectra $MSpin$, and of M Hopkins,
M Hovey \cite{HH}.

Let $KO_*(\cdot)$ be a periodic homological real $K$--theory, $KO$ be a
corresponding $\Omega$--spectrum. Also let $ko$ be the connected cover
of $KO$, and $ko\<2\>$ denote the $2$--connective cover of $ko$.  It is
convenient to identify the $2n$--fold connective covers of the spectrum
$KO$. Indeed, the $4k$--fold connective cover of $KO$ is
$\Sigma^{4k}ko$ (when $k$ is even), and the $(4k-2)$--fold connective
cover is $\Sigma^{4k-2} ko\<2\>$.  Let $ku$ be a connected cover of
the complex $K$--theory spectrum $K$. Let $\H(\Z/2)$ denote the
$\Z/2$--Eilenberg--MacLane spectrum. Recall that $ko$ and $ku$ are the
ring spectra with the coefficient rings:
\begin{equation}\label{ko_*}
\begin{array}{l}
ko_* \cong \Z[\eta, \omega, b]/(2\eta, \eta^3, \omega\eta, \omega^2-4b),
\ \ \ \deg\eta = 1, \ \deg \omega = 4, \deg b = 8; \!\!\!\!\!\!\!
\\
\\
ku_* \cong \Z[v ], \ \ \ \deg v = 2.
\end{array}
\end{equation}
Let $I = (i_1,\ldots,i_r)$ be a partition (possibly empty) of $n =
n(I)=\sum_{t=1}^r i_t$, $i_t > 0$. Each partition $I$ defines a map
$\pi^{I} \co MSpin \longrightarrow KO$ (which gives the
$KO$--characteristic class, see \cite{ABP}). If $I=\emptyset$ we denote
$\pi^{\emptyset}$ by $\pi^0$, which coincides with the
Atiyah--Bott--Shapiro orientation $\alpha \co MSpin \longrightarrow KO$.
\begin{Remark}
Let ${\cal P}$ be a set of all partitions, which is an abelian group.
We can make the set $\Z[{\cal P}]$ of linear combinations into a ring,
where multiplication of partitions is defined by set union, and then
to into a Hopf algebra with the diagonal $\displaystyle \Delta(I) =
\sum_{I_1+I_2=I} I_1\otimes I_2$.
\end{Remark}
Let $\mu \co MSpin\wedge MSpin \longrightarrow MSpin$, $\mu^{\prime} \co
KO\wedge KO \longrightarrow KO$ denote the ring spectra
multiplications. The Cartan formula says that
\begin{equation}\label{Cartan}
\begin{CD}
MSpin\wedge MSpin @>\mu>> MSpin \\
@V\sum(\pi^{I_1}\wedge\pi^{I_2})VV @V\pi^IVV \\
KO\wedge KO @>\mu^{\prime} >> KO
\end{CD} \ \ \ \mbox{or} \ \ \ \pi^{I_1} \mu = \sum_{I_1+I_2=I}\mu^{\prime}
(\pi^{I_1}\wedge\pi^{I_2}).
\end{equation}
\begin{Theorem}\label{ABP-Th}{\rm \cite{ABP}}
\begin{description}
\item{{\rm (1)}} Let $1\notin I$. \ Then if $n(I)$  is even, the
map $\pi^I \co MSpin \longrightarrow  KO$ lifts to a map $\bar{\pi}^I \co MSpin
\longrightarrow  \Sigma^{4n(I)}ko$.  If $n(I)$ is odd, the map $\pi^I$ lifts to a
map $\bar{\pi}^I \co MSpin \longrightarrow  \Sigma^{4n(I)-4}ko\<2\>$.
\item{{\rm (2)}} There exist a countable collection $z_k\in H^*(MSpin;\Z/2)$
such that the map 
$$
\begin{array}{l}
\displaystyle
\!\!\!\!\!\!\!\!\!\!\!\!\!\!\!\!
\prod_{1\notin
I}\!\bar{\pi}^I\!\!\times\!\prod_k \!z_k\co MSpin \longrightarrow \!\!\!\!\!\!\!\!\!\!
\prod_{\begin{array}{c}\scriptstyle1\notin I,\\ 
\scriptstyle n(I) \ {\rm even}\end{array}} 
\!\!\!\!\!\!\!\!\!\!\Sigma^{4n(I)}ko \times\!\!\!\!\!\!\!\!\!\!
\prod_{\begin{array}{c}\scriptstyle1\notin I,\\\scriptstyle n(I) \ 
{\rm odd}\end{array}} 
\!\!\!\!\!\!\!\!\!\!\Sigma^{4n(I)-4}ko\<2\> \times
\prod_{k}\Sigma^{\deg z_k} \H(\Z/2)
\end{array}
$$
is a 2--local homotopy equivalence.
\end{description}
\end{Theorem}
We use here the product symbol, however in the stable category of
spectra the product and the coproduct, ie the wedge, are the same.
We denote by $\rho^I$ the left inverses of the maps $\bar{\pi}^I$
(when $1\notin I$).
We denote also by $b$ an element in $\Omega^{Spin}_8$ which is the image
of the Bott element under the map $\rho^0$. The following Lemma due to
M Hovey and M Hopkins \cite{HH}. Since some fragments of its proof
will be used later, we provide an argument which essentially repeats
\cite{HH}.
\begin{Lemma}\label{b-loc}{\rm \cite[Lemma 1]{HH}} 
Let $I$ be a partition. Then $\pi^I(b)=0$ except for $\pi^0(b) =b$ and
possibly $\pi^{1}(b) \in KO_8$ and $\pi^{1,1}(b)\in KO_8$. The
elements $\pi^{1}(b)$, $\pi^{1,1}(b)$ are divisible by two in the
group $KO_8$. Further, the image of the Bott element $b$ is zero in
$MO_8$.
\end{Lemma}
\begin{proof} In the case $1\notin I$, $I\neq \emptyset$, the splitting
shows that $\pi^I(b) =0$. The map $\pi^I \co MSpin \longrightarrow KO$
(for any partitions $I$) may be lifted to the $4n(I)$--connective cover
of $KO$, as it is shown in \cite{Stong}. Let $S^0 \longrightarrow ko$
be a unit map, and $\rho^0 \co ko \longrightarrow MSpin$ be a left
inverse of $\bar{\pi}^0$. The composition
$$
\begin{CD}
S^0 \longrightarrow  ko @>\rho^0>> MSpin @>\pi^I>> KO
\end{CD}
$$
is null-homotopic for $I\neq \emptyset$. Let $\eta\in MSpin_1 =\Z/2$
be a generator. It is well-known that the image of the map $S^0
\longrightarrow MSpin$ on positive dimensional homotopy groups is
$\{b^n\eta, b^n\eta^2 \ | \ n\geq 0 \}$. It implies that
$\pi^I\rho^0(b^n\eta)=0$ and $\pi^I\rho^0(b^n\eta^2)=0$ for all
partitions $I\neq \emptyset$. Since the unit map $S^0\longrightarrow
MSpin$ is a map of ring spectra, we have
$\eta\cdot\pi^I\rho^0(b^n)=0$, so the elements $\pi^I\rho^0(b^n)$ are
even for all partitions $I\neq \emptyset$. In particular, $\pi^I(b)$
is even for all $I\neq \emptyset$.

Let $p_I$ be the Pontryagin class corresponding to a partition $I$.
Anderson, Brown and Peterson show that the Chern character $\mbox{ch
}(\pi^I(x)\otimes\C)= p_I(x) + \mbox{(higher terms)}$, for $x\in
\Omega_{Spin}^*(X)$. It implies that $p_I(b)$ are even elements for
all $I\neq \emptyset$. The Pontryagin classes $p_2$ and
$p_{1,1}=p_1^2$ determine the oriented cobordism ring $\Omega^{SO}_*$
in dimension 8, so the Bott element goes to an even element in
$\Omega^{SO}_8$ under the natural map $MSpin \longrightarrow MSO$.
Thus the composition $MSpin \longrightarrow MSO \longrightarrow MO$
takes the Bott element $b$ to zero.\hfill 
\end{proof}
We define the $K$--theory spectra with singularities 
$KO^{\Sigma_1}$, $KO^{\eta}$ and
$KO^{\Sigma_2}$, as the cofibers:
$$
\begin{array}{lcl}
\strut0{10}KO\wedge S^0 \stackrel{1\wedge 2}{\longrightarrow} KO\wedge S^0 
 \stackrel{\pi}{\longrightarrow} KO\wedge M(2) & =
& KO^{\Sigma_1} 
\\ 
\strut0{10}KO\wedge S^1 \stackrel{1\wedge \eta}{\longrightarrow} KO\wedge S^0 
\stackrel{\pi}{\longrightarrow}
KO\wedge \Sigma^{-2}\C\P^2& = & KO^{\eta}
\\
KO^{\eta}\wedge S^0 \stackrel{1\wedge 2}{\longrightarrow} 
KO^{\eta}\wedge S^0 \stackrel{\pi}{\longrightarrow}
KO^{\eta}\wedge M(2) & = & KO^{\Sigma_2}
\end{array}
$$
It is easy to derive (see, for example, \cite{Mah2}) 
the following statement.
\begin{Corollary}\label{KO-2}
The spectrum $KO^{\eta}$ is homotopy equivalent (as a ring spectrum)
to the spectrum $K$, classifying the complex $K$--theory, and the
spectrum $KO^{\Sigma_2}$ is homotopy equivalent (as a ring spectrum)
to the spectrum $K(1)$ classifying the first Morava $K$--theory.
\end{Corollary}
We introduce also the notation:
$$
\begin{array}{ll}
\strut0{10}\!ko^{\Sigma_1} = ko\wedge M(2), 
& \!\!\!ko\<2\>^{\Sigma_1} = ko\<2\> \wedge M(2),
\H(\Z/2)^{\Sigma_1} = \H(\Z/2)\wedge M(2);
\\
\strut0{10}\!ko^{\eta}\! = ko\wedge \Sigma^{-2}\C\P^2, \!\!\!\!\!\!\!\!\!
& \!\!\!ko\<2\>^{\eta} \!= \!ko\<2\> \!\wedge \!\Sigma^{-2}\C\P^2, \ 
\H(\Z/2)^{\eta}\! = \!\H(\Z/2)\wedge \Sigma^{-2}\C\P^2;
\\
\!ko^{\Sigma_2} = ko^{\eta}\wedge M(2), 
& \!\!\!ko\<2\>^{\Sigma_2} = ko\<2\>^{\eta} \wedge M(2),
\H(\Z/2)^{\Sigma_2} = \H(\Z/2)^{\eta}\wedge M(2).
\end{array}
$$
Let $I$ be a partition as above. The $KO$--characteristic numbers
$$
\pi^I \co MSpin \longrightarrow KO
$$ 
which are lifted to the connective
cover $ko\<4n(I)\>$ give the characteristic numbers
$$
\begin{array}{llll}
\strut0{10}\pi^I_{\Sigma_1}\!\!\!&=&\!\!\!\pi^I\wedge 1 \co 
MSpin^{\Sigma_1} = MSpin\wedge M(2) \longrightarrow 
KO\wedge M(2)\!\!\! & = KO^{\Sigma_1},
\\
\strut0{10}\pi^I_{\eta}\!\!\!&=&\!\!\!\pi^I\wedge 1 \co 
MSpin^{\eta} = MSpin\wedge \Sigma^{-2}\C\P^2
\longrightarrow 
KO\wedge \Sigma^{-2}\C\P^2\!\!\! &= KO^{\eta},
\\
\pi^I_{\Sigma_2}\!\!\!&=&\!\!\!\pi^I\wedge 1 \co 
MSpin^{\Sigma_1} = MSpin\wedge Y \longrightarrow 
KO\wedge Y\!\!\! &= KO^{\Sigma_2},
\end{array}
$$ 
together with the lifts to the corresponding connective covers:
$$
\begin{array}{lllll}
\strut0{10}\bar\pi^I_{\Sigma_1} \!\!\!& = & \!\!\!\bar\pi^I\wedge 1 \co 
&MSpin^{\Sigma_1} \longrightarrow 
ko\<4n(I)\> \wedge M(2)& = ko\<4n(I)\>^{\Sigma_1}
\\
\strut0{10}\bar\pi^I_{\eta}\!\!\!& = &\!\!\!\bar\pi^I\wedge 1 \co 
&MSpin^{\eta} 
\longrightarrow 
ko\<4n(I)\>\wedge \Sigma^{-2}\C\P^2& = ko\<4n(I)\>^{\eta}
\\
\bar\pi^I_{\Sigma_2}\!\!\!& = & \!\!\!\bar\pi^I\wedge 1 \co 
&MSpin^{\Sigma_1} \longrightarrow 
ko\<4n(I)\>\wedge Y & = ko\<4n(I)\>^{\Sigma_2}
\end{array}
$$ 
Now we would like to identify the spectra $ko^{\Sigma}$,
$ko\<4n(I)\>^{\Sigma}$ for $\Sigma=\Sigma_2$ or $\eta$ for those
partitions $I$, $1\notin I$. It is enough to determine a homotopy type
of the spectra $ko^{\Sigma}$ and $ko\<2\>^{\Sigma}$.

Let $\A(1)$ be a subalgebra of the Steenrod algebra $\A_2$ generated
by $1, Sq^1, Sq^2$. The cohomology $H^*(ko)$ as a module over Steenrod
algebra is $ H^*(ko) \cong \A_2\otimes_{\A(1)}\Z/2 $. The K\"unneth
homomorphism
$$
H^*(ko\wedge X)  \cong (\A_2\otimes_{\A(1)}\Z/2) \otimes H^*(X) \cong 
\A_2\otimes_{\A(1)} H^*(X) 
$$ 
and the ring change formula $ \Hom_{\A_2}(\A_2\otimes_{\A(1)} M,
N)\cong \Hom_{\A(1)}(M,N) $ turn the ordinary mod 2 Adams spectral
sequence into the one with the $E_2$--term
$$
\Ext_{\A(1)}^{s,t}(H^*(X),\Z/2) \Longrightarrow ko_{t-s}(X).
$$
Here we use regular conventions to draw the cell-diagrams for
the spectra in question. Recall that

$$
H^*(ko) =  \A_2\otimes_{\A(1)} 
\mbox{
\begin{picture}(2,20)
\put(0,0){\circle*{3}}
\end{picture}} \ \ \ \mbox{and}\ \ \ 
H^*(ko\<2\>)= \A_2\otimes_{\A(1)} 
\mbox{
\begin{picture}(2,20)
\put(0,0){\circle*{3}}
\put(0,20){\circle*{3}}
\put(0,10){\circle*{3}}
\put(0,30){\circle*{3}}
\put(0,40){\circle*{3}}
\put(0,0){\line(0,1){10}}
\put(0,30){\line(0,1){10}}
\put(0,10){\oval(10,20)[br]}
\put(0,10){\oval(10,20)[tr]}
\put(0,30){\oval(10,20)[br]}
\put(0,30){\oval(10,20)[tr]}
\put(0,20){\oval(10,20)[bl]}
\put(0,20){\oval(10,20)[tl]}
\end{picture}} \ \ \ \mbox{(the joker).}
$$
Let $k(1)$ be a connected cover of the first Morava $k$--theory
spectrum $K(1)$ with the coefficient ring $k(1)_*\cong
\Z/2[v_1]$. Here is the result for the spectra $ko^{\eta}$,
$ko^{\Sigma_2}$:
\begin{Lemma}\label{ko-Sigma}
There are the following homotopy equivalences
\begin{equation}\label{eq5a}
\begin{array}{l}
ko^{\eta} \ \ \cong \ \ ku, \ \ \ \ \ 
ko^{\Sigma_2} \ \cong \ \ k(1)
\end{array}
\end{equation}
\end{Lemma}
The following result one can prove by an easy computation:
\begin{Lemma} There are isomorphisms of the following $\A(1)$--modules:
\begin{equation}\label{eq6}
\mbox{
\begin{picture}(300,60)
\put(50,0){\circle*{3}}
\put(50,20){\circle*{3}}
\put(50,10){\circle*{3}}
\put(50,30){\circle*{3}}
\put(50,40){\circle*{3}}
\put(50,0){\line(0,1){10}}
\put(50,30){\line(0,1){10}}
\put(50,10){\oval(10,20)[br]}
\put(50,10){\oval(10,20)[tr]}
\put(50,30){\oval(10,20)[br]}
\put(50,30){\oval(10,20)[tr]}
\put(50,20){\oval(10,20)[bl]}
\put(50,20){\oval(10,20)[tl]}
\put(40,-3){{\small $a$}}
\put(40,7){{\small $b$}}
\put(57,17){{\small $c$}}
\put(40,27){{\small $d$}}
\put(40,37){{\small $e$}}

\put(90,0){\circle*{3}}
\put(90,20){\circle*{3}}
\put(90,10){\oval(10,20)[br]}
\put(90,10){\oval(10,20)[tr]}
\put(80,-3){{\small $\alpha$}}
\put(80,17){{\small $\beta$}}

\put(65,10){{\large $\otimes$}}
\put(105,10){{\large $=$}}

\put(130,10){\circle*{3}}
\put(130,0){\line(0,1){9}}
\put(130,0){\circle*{3}}
\put(130,20){\oval(20,40)[br]}
\put(130,30){\oval(20,40)[bl]}
\put(140,20){\circle*{3}}
\put(140,20){\line(0,1){10}}
\put(120,30){\circle*{3}}
\put(120,30){\line(0,1){10}}
\put(120,20){\line(0,1){10}}
\put(140,30){\circle*{3}}
\put(120,40){\circle*{3}}
\put(130,30){\oval(20,40)[tr]}
\put(140,20){\line(-1,1){20}}
\put(130,40){\oval(20,40)[tl]}
\put(130,50){\circle*{3}}
\put(130,50){\line(0,1){9}}
\put(130,60){\circle*{3}}

\put(117,-3){{\small $a\alpha$}}
\put(150,10){{\large $\oplus$}}

\put(170,20){\circle*{3}}
\put(170,40){\circle*{3}}
\put(170,30){\oval(10,20)[br]}
\put(170,30){\oval(10,20)[tr]}
\put(170,13){{\small $c\alpha$}}
\end{picture}}
\end{equation}
\begin{equation}\label{eq7}
\mbox{
\begin{picture}(300,60)
\put(50,0){\circle*{3}}
\put(50,20){\circle*{3}}
\put(50,10){\circle*{3}}
\put(50,30){\circle*{3}}
\put(50,40){\circle*{3}}
\put(50,0){\line(0,1){10}}
\put(50,30){\line(0,1){10}}
\put(50,10){\oval(10,20)[br]}
\put(50,10){\oval(10,20)[tr]}
\put(50,30){\oval(10,20)[br]}
\put(50,30){\oval(10,20)[tr]}
\put(50,20){\oval(10,20)[bl]}
\put(50,20){\oval(10,20)[tl]}
\put(40,-3){{\small $a$}}
\put(40,7){{\small $b$}}
\put(57,17){{\small $c$}}
\put(40,27){{\small $d$}}
\put(40,37){{\small $e$}}

\put(90,0){\circle*{3}}
\put(90,20){\circle*{3}}
\put(90,10){\circle*{3}}
\put(90,30){\circle*{3}}
\put(90,0){\line(0,1){10}}
\put(90,20){\line(0,1){10}}
\put(90,10){\oval(10,20)[br]}
\put(90,10){\oval(10,20)[tr]}
\put(90,20){\oval(10,20)[bl]}
\put(90,20){\oval(10,20)[tl]}
\put(78,-3){{\small $\alpha$}}
\put(78,7){{\small $\beta$}}
\put(95,17){{\small $\gamma$}}
\put(95,27){{\small $\delta$}}

\put(65,10){{\large $\otimes$}}
\put(105,10){{\large $=$}}

\put(130,10){\circle*{3}}
\put(130,0){\line(0,1){9}}
\put(130,0){\circle*{3}}
\put(130,20){\oval(20,40)[br]}
\put(130,30){\oval(20,40)[bl]}
\put(140,20){\circle*{3}}
\put(140,20){\line(0,1){10}}
\put(120,30){\circle*{3}}
\put(120,30){\line(0,1){10}}
\put(120,20){\line(0,1){10}}
\put(140,30){\circle*{3}}
\put(120,40){\circle*{3}}
\put(130,30){\oval(20,40)[tr]}
\put(140,20){\line(-1,1){20}}
\put(130,40){\oval(20,40)[tl]}
\put(130,50){\circle*{3}}
\put(130,50){\line(0,1){9}}
\put(130,60){\circle*{3}}

\put(117,-3){{\small $a\alpha$}}

\put(150,10){{\large $\oplus$}}

\put(170,20){\circle*{3}}
\put(170,10){\line(0,1){9}}
\put(170,10){\circle*{3}}
\put(170,30){\oval(20,40)[br]}
\put(170,40){\oval(20,40)[bl]}
\put(180,30){\circle*{3}}
\put(180,30){\line(0,1){10}}
\put(160,40){\circle*{3}}
\put(160,40){\line(0,1){10}}
\put(160,30){\line(0,1){10}}
\put(180,40){\circle*{3}}
\put(160,50){\circle*{3}}
\put(170,40){\oval(20,40)[tr]}
\put(180,30){\line(-1,1){20}}
\put(170,50){\oval(20,40)[tl]}
\put(170,60){\circle*{3}}
\put(170,60){\line(0,1){9}}
\put(170,70){\circle*{3}}

\put(165,3){{\small $b\alpha$}}

\put(190,10){{\large $\oplus$}}
\put(205,13){{\small $c\alpha$}}

\put(210,20){\circle*{3}}
\put(210,40){\circle*{3}}
\put(210,30){\circle*{3}}
\put(210,50){\circle*{3}}
\put(210,20){\line(0,1){10}}
\put(210,40){\line(0,1){10}}
\put(210,30){\oval(10,20)[br]}
\put(210,30){\oval(10,20)[tr]}
\put(210,40){\oval(10,20)[bl]}
\put(210,40){\oval(10,20)[tl]}
\end{picture}}
\end{equation}
\end{Lemma}
Using the Adams spectral sequence for the spectra $ko\<2\>^{\eta}$ and
$ko\<2\>^{\Sigma_2}$, one obtains the following result:
\begin{Lemma}
There are the following homotopy equivalences
\begin{equation}\label{eq8}
\begin{array}{ll}
\strut08 ko\<2\>^{\eta} \ \ \cong & \H(\Z/2)\vee \Sigma^2 ku,
\\
ko\<2\>^{\Sigma_2} \ \cong & \H(\Z/2) \vee \Sigma \H(\Z/2)\vee  
\Sigma^2 k(1). 
\end{array}
\end{equation}
\end{Lemma}
It is convenient to denote:
$$
\widehat{ko}= \!\!\!\!\!\!\!\!\!\!\!\!\!\!
\prod_{\begin{array}{c}\scriptstyle1\notin I,\\ 
\scriptstyle n(I)\neq 0, \ {\rm even}\end{array}} 
\!\!\!\!\!\!\!\!\!\!\!\!\!\!\Sigma^{4n(I)}ko \times\!\!\!\!\!\!\!\!\!\!
\prod_{\begin{array}{c}\scriptstyle1\notin I,\\\scriptstyle n(I) \ 
{\rm odd}\end{array}} 
\!\!\!\!\!\!\!\!\!\!\Sigma^{4n(I)-4}ko\<2\>, \ \ \ \mbox{and} \ \ \ 
\widehat{\H(\Z/2)}=\prod_{k}\Sigma^{\deg z_k} \H(\Z/2).
$$
The spectra $\widehat{ko}^{_{\Sigma}}$ are defined similarly for
$\Sigma= \Sigma_1$, $\Sigma_2$, $\Sigma_3$ or $\eta$.  Theorem
\ref{ABP-Th} implies the following result:
\begin{Corollary}\label{ABP-Sigma}
There is the following homotopy equivalence 
of 2--local spectra:
$$
F^{\Sigma}\co MSpin^{{\Sigma}} \longrightarrow  
ko^{\Sigma}\vee \widehat{ko}^{{\Sigma}}\vee
\widehat{\H(\Z/2)}^{{\Sigma}}, \ \ \ 
\mbox{where $\Sigma= \Sigma_1$, $\Sigma_2$, or $\eta$.}
$$
\end{Corollary}
\begin{Remark}
The coefficient groups of the $K$--theories $KO^{\Sigma}$ are well-known
in homotopy theory. We give the table of the groups $KO^{\Sigma_1}_n=KO_n(pt,\Z/2)$
for convenience:

\begin{small}

\begin{tabular}{|l|c|c|c|c|c|c|c|c|c|c|}\hline
\strut{10}3  &   $0$  &   $1$  & $2$ & $3$& $4$  & $5$ &$6$ & $7$ & $8$ &
\\ \hline
\strut{12}5$KO^{\Sigma_1}_n=KO_n(pt,\Z/2)$&$\Z/2$&$\Z/2$&$\Z/4$&$\Z/2$&$\Z/2$&$0$&$0$&$0$&$\Z/2$ &$\cdots$
\\ \hline
\end{tabular}
\end{small}

We emphasize that $KO_{8k+2}(pt,\Z/2)\cong \Z/4$.
\end{Remark}
\begin{Remark} We notice that there is a natural transformation 
$$ r :
\Omega^{Spin,\eta}_*(\cdot) \longrightarrow  \Omega^{Spin^c}_*(\cdot).  
$$ 
Indeed,
let $M$ be an $\eta$--manifold, ie, $\p M \cong \beta_2 M \times P_2$,
where $P_2=S^1$ with nontrivial $Spin$ structure. Then $P_2$ is a
boundary as a $Spin^c$--manifold, even more, $P_2=\p D^2$.  Then the
correspondence
$$
(M, \p M = \beta_2 M\times P_2) \mapsto (N = M\cup -\beta_2 M \times D^2)
$$
determines the transformation $r$. In particular, $r$ gives a map of
classifying spectra: $r \co MSpin^{\eta} \longrightarrow MSpin^c$. It is
easy to see that there is a commutative diagram
$$
\begin{diagram}
\setlength{\dgARROWLENGTH}{1.2em}
	\node{MSpin^{\eta}}
		\arrow[4]{e,t}{r}
		\arrow{s,l}{\cong}
	\node[4]{MSpin^c}
		\arrow{s,l}{\cong}
\\
	\node{MSpin \wedge \Sigma^{-2}\C\P^2}
		\arrow[4]{e,t}{Id \wedge \Sigma^{-2}j}
	\node[4]{MSpin\wedge \Sigma^{-2}\C\P^{\infty}}
\end{diagram}
$$
where $j \co \C\P^2 \longrightarrow \C\P^{\infty}$ is the standard
embedding.  There are simple geometric reasons which imply that the
transformation $r$ is not multiplicative.  In fact, it is very similar
to the transformation $ \Omega^{SU,\eta}_*(\cdot) \longrightarrow
\Omega^{U}_*(\cdot) $, where $\Omega^{SU,\eta}_*(\cdot)$ is the
$SU$--cobordism theory with $\eta$--singularities. The cobordism theory
$\Omega^{SU,\eta}_*(\cdot)$ may be easily identified with the
Conner--Floyd theory $W(\C,2)_*(\cdot)$, see {\rm \cite{Mir}}.
\end{Remark}
\section{The spectrum $MSpin^{\Sigma_3}$}\label{s6}
Let $A \co \Sigma^8 M(2) \longrightarrow  M(2)$ be the Adams map.  
Let $V(1)$ be a cofiber:
$$
\begin{CD}
\Sigma^8 M(2) @>A>> M(2) @>p>> V(1).
\end{CD}
$$
The objective of this section is to prove the following result.
\begin{Theorem}\label{Sigma3}
There is a homotopy equivalence of spectra localized at 2:
\begin{equation}\label{eq10}
MSpin^{\Sigma_3}\cong MSpin \wedge
\Sigma^{-2}\C\P^2\wedge V(1).
\end{equation}
\end{Theorem}
\begin{proof} Recall that the Adams map $A$ induces a multiplication
by the Bott element in $KO_*$ and connected covers $ko_*$
and $ko\<2\>_*$. Let, as above,
$Y=\Sigma^{-2}\C\P^{2}\wedge M(2)$. We apply the Cartan formula
(\ref{Cartan})
$$
\begin{CD}
MSpin\wedge MSpin \wedge Y @>\mu\wedge 1>> MSpin \wedge Y \\
@V\sum(\pi^{I_1}\wedge\pi^{I_2})\wedge 1VV @V\pi^IVV \\
KO\wedge KO \wedge Y @>\mu^{\prime}\wedge 1 >> KO\wedge Y
\end{CD} 
$$
to obtain the formula:
\begin{equation}\label{Cartan1}
\begin{CD}
MSpin \wedge MSpin^{\Sigma_2} @>\mu_{\Sigma_2} >> MSpin^{\Sigma_2} \\
@V\sum(\pi^{I_1}\wedge\pi^{I_2}_{\Sigma_2})VV @V\pi^I_{\Sigma_2}VV \\
KO \wedge KO^{\Sigma_2} @>\mu^{\prime}_{\Sigma_2} >> KO^{\Sigma_2}
\end{CD} \!\mbox{or} \ \
\mu_{\Sigma_2}\pi^{I}_{\Sigma_2} = \sum_{I_1+I_2=I}\mu^{\prime}_{\Sigma_2}
(\pi^{I_1}\wedge\pi^{I_2}_{\Sigma_2}).
\end{equation}
Let $X$ be a space, $x\in \Omega^{Spin}_*(X)$, $b\in \Omega^{Spin}_8$
be the Bott element. Then $\mu_{\Sigma_2}(b,x)=b\cdot x\in
\Omega^{Spin}_{*+8}(X)$.
\begin{Lemma}\label{bott2}
The $KO^{\Sigma_2}$--characteristic numbers $\pi^{I}_{\Sigma_2} \co
MSpin^{\Sigma_2} \longrightarrow KO^{\Sigma_2}$ commutes with a
multiplication by the Bott element, ie, $ \pi^{I}_{\Sigma_2}(b\cdot
x) = b\cdot \pi^{I}_{\Sigma_2}(x) $.
\end{Lemma}
\begin{proof} The Cartan formula (\ref{Cartan1}) and Lemma
\ref{b-loc} gives:
$$
\begin{array}{l}
\displaystyle
\strut0{18}\pi^{I}_{\Sigma_2}(b\cdot x) = 
\sum_{I_1+I_2=I}(\pi^{I_1}(b)\pi^{I_2}_{\Sigma_2}(x) 
= b \cdot \pi^{I}_{\Sigma_2}(x) + \pi^{(1)}(b)y +
\pi^{(1,1)}(b)z = 
\\
b\cdot \pi^{I}_{\Sigma_2}(x) + (2c)\cdot y + (2d)\cdot z = 
b\cdot \pi^{I}_{\Sigma_2}(x) + c\cdot (2 y) + d\cdot (2z) =
b\cdot \pi^{I}_{\Sigma_2}(x).
\end{array}
$$
Here $\pi^{(1)}(b)=2c$, $y,z \in KO_*^{\Sigma_2}(X)$, and
$\pi^{(1,1)}(b)=2d$ by Lemma \ref{b-loc}. We note that $2 y =0$ and
$2z=0$ since the cobordism theory $\Omega^{Spin,\Sigma_2}_*(\cdot)$
has an admissible product structure by Theorem \ref{product}. \hfill
\end{proof}
Let $I$ be a partition, and $1\notin I$. The map $\pi^{I}_{\Sigma_2} \co
MSpin_{\Sigma_2} \longrightarrow KO_{\Sigma_2}$ lifts to connective
cover: $\bar\pi^{I}_{\Sigma_2} \co MSpin^{\Sigma_2} \longrightarrow
ko\<4n(I)\>^{\Sigma_2}$.  Let $S^8 \stackrel{b}{\longrightarrow}
MSpin$ be a map representing the Bott element $b\in
\Omega^{Spin}_8$. We denote by $\cdot b$ the composition
$$
\begin{CD}
\Sigma^8 MSpin = S^8\wedge MSpin^{\Sigma_2} 
@>b\wedge 1>> MSpin \wedge MSpin^{\Sigma_2}
@>\mu^{\prime\prime}>> MSpin^{\Sigma_2} .
\end{CD}
$$
Note that the diagram
\begin{equation}\label{eq11a}
\begin{diagram}
\setlength{\dgARROWLENGTH}{0.8em}
	\node{\Sigma^8 MSpin^{\Sigma_2}}
		\arrow[6]{e,t}{\cdot b}
		\arrow{s,l}{\cong}
	\node[6]{MSpin^{\Sigma_2}}
		\arrow{s,l}{\cong}
\\
	\node{\Sigma^8 MSpin\wedge  M(2)\wedge \Sigma^{-2}\C\P^2}
		\arrow[6]{e,t}{\cdot b\wedge 1 \wedge 1}
	\node[6]{MSpin\wedge  M(2)\wedge \Sigma^{-2}\C\P^2}
\end{diagram}
\end{equation}
commutes since $MSpin^{\Sigma_2}$ is a module over $MSpin$.
\begin{Lemma}\label{Lbott}
Let $I$ be a partition, so that $1\notin I$.
The following diagrams commute:
\begin{equation}\label{eq12}
\!\!\!\!\!\!\!\!\!\!\!\!\!\!\!\!\!\!
\begin{diagram}
\setlength{\dgARROWLENGTH}{0.6em}
	\node{\Sigma^8 MSpin^{\Sigma_2}}
		\arrow[6]{e,t}{\cdot b}
		\arrow{s,r}{\Sigma^8 \pi^I_{\Sigma_2}}
	\node[6]{MSpin^{\Sigma_2}}
		\arrow{s,r}{\pi^I_{\Sigma_2}}
\\
	\node{\Sigma^8 ko\<n(I)\>^{\Sigma_2}}
		\arrow{s,r}{\cong}
	\node[6]{ko\<n(I)\>^{\Sigma_2}}
		\arrow{s,r}{\cong}
\\
	\node{ko\<n(I)\>\wedge \Sigma^8 M(2)\wedge \Sigma^{-2}\C\P^2}
		\arrow[6]{e,t}{1\wedge A \wedge 1}
	\node[6]{ko\<n(I)\>\wedge M(2)\wedge \Sigma^{-2}\C\P^2\!\!\!\!\!\!\!}
\\
	\node{\Sigma^8 MSpin^{\Sigma_2}}
		\arrow[6]{e,t}{\cdot b}
		\arrow{s,r}{\Sigma^8 z_k^{\Sigma_2}}
	\node[6]{MSpin^{\Sigma_2}}
		\arrow{s,r}{\Sigma^8 z_k^{\Sigma_2}}
\\
	\node{\Sigma^{8+\deg z_k}\H(\Z/2)\wedge M(2)\wedge \Sigma^{-2}\C\P^2}
		\arrow{s,r}{\cong}
	\node[6]{\Sigma^{\deg z_k}\!\H(\Z/2)\!\wedge \!M(2)\!\wedge\!\Sigma^{-2}\!\C\P^2\!\!\!\!\!\!\!\!\!\!\!\!\!\!\!\!\!\!\!\!\!}
		\arrow{s,r}{\cong}
\\
	\node{\Sigma^{\deg z_k}\!\H(\Z/2)\!\wedge \!\Sigma^8\!
M(2)\!\wedge \!\Sigma^{-2}\!\C\P^2}
		\arrow[6]{e,t}{1\wedge A\wedge 1}
	\node[6]{\Sigma^{\deg z_k}\!\H(\Z/2)\wedge \!M(2)\!\wedge\!
\Sigma^{-2}\!\C\P^2\!\!\!\!\!\!\!\!\!\!\!\!\!\!\!\!\!\!\!\!}
\end{diagram}\!\!\!\!\!\!\!\!\!\!
\end{equation}
\end{Lemma}
\begin{proof} A commutativity of the first diagram follows from Lemma
\ref{bott2} and the diagram (\ref{eq11a}). Recall that a projection
of the Bott element into the homotopy group of $\Sigma^{\deg
z_k}\H(\Z/2)$ is zero.  Let $X$ be a finite spectrum.
The map 
$$
\begin{array}{l}
\strut0{10}1\wedge A \wedge 1\wedge 1\co \Sigma^{\deg z_k}\H(\Z/2)\wedge \Sigma^8
M(2) \wedge \Sigma^{-2} \C\P^2 \wedge X 
\longrightarrow \mbox{ \hspace*{10mm} \ } 
\\
\mbox{ \hspace*{30mm} \ } \ \ \ \ \ \ \ \ \ \ \ \ 
\Sigma^{\deg z_k} \H(\Z/2)\wedge  M(2) \wedge
\Sigma^{-2} \C\P^2 \wedge X 
\end{array}
$$
in homotopy coincides with the homomorphism in mod 2 homology groups
$$
\begin{CD}
\Sigma^{\deg z_k}H_*(\Sigma^8\!
M(2)\!\wedge \!\Sigma^{-2}\!\C\P^2 \wedge X) 
@>A_* \otimes 1\otimes  1>>
\Sigma^{\deg z_k}H_*(M(2)\!\wedge \!\Sigma^{-2}\!\C\P^2 \wedge X)
\end{CD}
$$
and is trivial for any space $X$ since $A$ has the Adams filtration 4.
It implies that $1\wedge A \wedge 1$ is a trivial map. A commutativity
of (\ref{eq12}) now follows. \hfill
\end{proof}
To complete the proof of Theorem \ref{Sigma3} we notice that
Lemmas \ref{Lbott} and \ref{bott2} give the commutative 
diagram
$$
\begin{diagram}
\setlength{\dgARROWLENGTH}{1.2em}
	\node{\Sigma^8 MSpin^{\Sigma_2}}
		\arrow[3]{e,t}{\cdot b}
		\arrow{s,r}{\Sigma^8 F^{\Sigma_2}}
	\node[3]{MSpin^{\Sigma_2}}
		\arrow[3]{e,t}{\pi_3}
		\arrow{s,r}{F^{\Sigma_2}}
	\node[3]{MSpin^{\Sigma_3}}
		\arrow{s,r}{F^{\Sigma_3}}
\\
	\node{\Sigma^8 MSpin^{\eta}\wedge \Sigma^8 M(2)}
		\arrow[3]{e,t}{1\wedge A}
	\node[3]{MSpin^{\eta}\wedge M(2)}
		\arrow[3]{e,t}{p}
	\node[3]{MSpin^{\eta}\wedge V(1)}
\end{diagram}
$$
where the map $F^{\Sigma_3}$ exists since the both rows are
cofibrations. The five-lemma implies that $F^{\Sigma_3}$
is a homotopy equivalence. \hfill 
\end{proof}
\begin{Corollary}\label{KO-3}
	The spectrum $KO^{\Sigma_3}= KO\wedge \Sigma^{-2}\C\P^2 \wedge
	V(1)$ is a contractible spectrum.
\end{Corollary}
\begin{Remark}
	The connective spectrum $ko^{\Sigma_3}$ is of some
	interest. It is certainly not contractible, and it is very
	easy to see that
$$
ko^{\Sigma_3}_j =\{ \begin{array}{cl} 
\Z/2 & \mbox{if $j =0, 2,4,6$,}
\\
0  & \mbox{otherwise,}
\end{array}\right.
$$
	and the Postnikov tower of $ko^{\Sigma_3}$ has the operation
	$Q_1$ as its $k$--invariants.
\end{Remark}
The technique we used above may be applied to prove the following result:
\begin{Corollary}\label{MSpin-eta}
	There is such admissible product structure $\mu^{(2)}$ of the
	spectrum $MSpin^{\Sigma_2}$, so that the map
	$\bar\pi^0_{\Sigma_2} \co MSpin^{\Sigma_2} \longrightarrow
	ko^{\Sigma_2}=k(1)$ is a ring spectra map, moreover, there is
	an inverse ring spectra map $\rho^0_{\Sigma_2} \co ko^{\Sigma_2}
	\longrightarrow MSpin^{\Sigma_2}$. In other words,
	$ko^{\Sigma_2}$ splits off of the spectrum $MSpin^{\Sigma_2}$
	as a ring spectrum.
\end{Corollary}
\section{Surgery Lemma for $\Sigma$--manifolds}\label{s7}
{\bf \ref{s7}.1\qua A Riemannian metric on a $\Sigma$--manifold}\qua  Here we
describe what do we mean by a Riemannian metric on manifold with
singularities. 
We consider the case when a manifold has of at most three
singularities, $\Sigma_3=(P_1,P_2,P_3)$.  We denote $\Sigma_1=(P_1)$,
$\Sigma_2=(P_1,P_2)$.  We assume that there are given Riemannian
metrics $g_{P_i}$ on the manifolds $P_i$, $i=1,2,3$. As we mentioned
earlier, the metrics $g_{P_i}$ are not assumed to be psc-metrics.

If $M$ is a $\Sigma_3$--manifold, we assume that it is given a
decomposition of the boundary $\p M$:
$$
\begin{array}{l}
\strut08
\p M = \(\beta_1 M \times P_1 
\cup \beta_2 M \times P_2  \cup \beta_2 M \times P_2\)
\cup \ \beta_{123} M \times P_1 \times P_3 \times P_2 \times D^2
\\
\cup \ \(\beta_{12} M \times P_1 \times P_2 \times I_{12}
\cup \beta_{23} M \times P_2 \times P_3 \times I_{23}
 \cup \beta_{13} M \times P_1 \times P_3 \times I_{13}\)
\end{array}
$$
glued together as it is shown on Figure 7 (a). We start with a Riemannian
metric $g_{123}$ on the manifold $\beta_{123} M$. We
assume that the manifold
$$
\beta_{123}M \times P_1 \times P_2 \times P_3 \times D^2
$$
has product metric $g_{123}\times g_{P_1}\times g_{P_2}\times
g_{P_3}\times g_0$, where $g_0$ is the standard flat metric on the
disk $D^2$. 

\noindent
Besides, we assume that the manifold $\beta_{123}M \times
P_1 \times P_2 \times P_3 $, being common boundary of the manifolds
$$
\beta_{12} M \times P_1 \times P_2 ,  \ \ \ 
\beta_{13} M \times P_1 \times P_3 , \ \  \ \mbox{and} \ \ \
\beta_{23} M \times P_2 \times P_3 ,
$$
is embedded together with the colors (see Figure 7 (a)):
$$
\begin{array}{l}
\strut0{10}\beta_{123}M \times P_1 \times P_2 \times P_3 \times I^{\prime}_{12}
\subset \beta_{12} M \times P_1 \times P_2 ,
\\
\strut0{10}\beta_{123}M \times P_1 \times P_2 \times P_3 \times I^{\prime}_{13}
\subset \beta_{13} M \times P_1 \times P_3,
\\
\beta_{123}M \times P_1 \times P_2 \times P_3 \times I^{\prime}_{23}
\subset \beta_{23} M \times P_2 \times P_3.
\end{array}
$$
Here $I^{\prime}_{ij}$ are the intervals embedded into the flat disk
$D^2$ as it is shown on Figure 7 (b).
\vspace{2mm}

\noindent
\hspace*{12mm}
{\PSbox{botv-6.ps}{5cm}{60mm}
\begin{picture}(0,0)
\put(-45,100){{\small $\beta_{12}M\!\times \!P_1\!\times \!P_2\!\times 
\!I_{12}$}}
\put(-160,-10){{\small $\beta_{13}M\!\times \!P_1\!\times \!P_3\!\times 
\!I_{13}$}}
\put(-173,160){{\small $\beta_{23}M\!\times \!P_2\!\times \!P_3\!\times \!
I_{23}$}}
\put(-75,165){{\small $\beta_{123}M\!\times \!P_1\!\times
\!P_2\!\times \!P_3\!\times \! D^2$}} 
\put(-55,35){{\small $\beta_1 M\!\times \! P_1$}}
\put(-63,125){{\small $\beta_2 M\!\times \! P_2$}} 
\put(-140,75){{\small $\beta_3 M\!\times \! P_3$}}
\put(-180,75){$\p M =$}
\put(80,-30){{\small  (b)\qua Flat disk $D^2$}}
\put(-160,-30){{\small (a)\qua The decomposition of $\p M$}}
\put(150,80){{\small $I_{12}^{\prime}$}}
\put(100,100){{\small $I_{23}^{\prime}$}}
\put(100,40){{\small $I_{13}^{\prime}$}}
\end{picture}}\hspace*{13mm}\PSbox{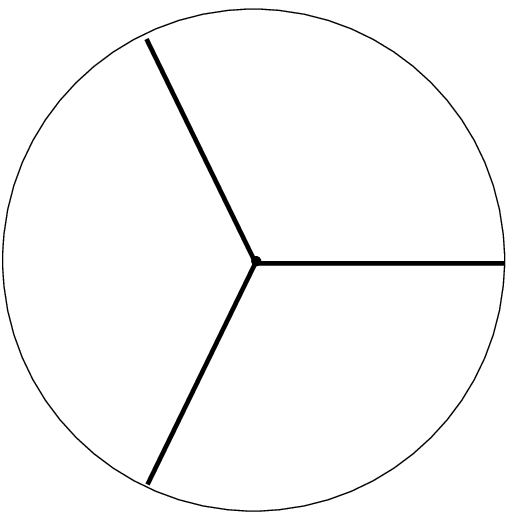}{5cm}{60mm}
\vspace*{10mm}

\centerline{{\small Figure 7}}

Let $g_{ij}$ be metrics on the manifolds $\beta_{ij} M$. We assume
that the product metric $$
g_{ij}\times g_{P_i}\times g_{P_j} \ \ \ \mbox{on the manifold}
\ \ \ \beta_{ij} M \times P_i \times P_j
$$ 
coincides with the product metric on the color $\beta_{123}M \times
P_1 \times P_2 \times P_3 \times I^{\prime}_{ij}$ near its boundary.
Finally if $g_i$ is a metric in $\beta_i M$ ($i=1,2,3$), then we
assume that
the product metric $g_i\times g_{P_i}$ on $\beta_i M \times P_i$
coincides with the above product metrics on the manifold $\beta_{ij} M
\times g_{P_i}\times g_{P_j}\times I_{ij}$. Furthermore, the product
metric $g_i\times g_{P_i}$ on $\beta_i M \times P_i$ restricted on the
manifold
$$
(\beta_{123}M \times P_1 \times P_2 \times P_3 \times D^2 ) \cap
(\beta_i M \times P_i)
$$
coincides with the product metric $g_{123}\times g_{P_1}\times
g_{P_2}\times g_{P_3}\times g_0$. Finally the metric $g$ on the manifold
$M$ is assumed to be product metric near the boundary $\p M$. 
Let $M$, as above be a $\Sigma$--manifold with the same singularities
$\Sigma=(P_1,P_2,P_3)$. We say that a {\sl metric $g$ on $M$ is of
positive scalar curvature}, if, besides the above conditions, the
metrics $g$ on $M$, $g_{i}$ on $\beta_i M$, $g_{ij}$ on $\beta_{ij}M$,
and $g_{123}$ on $\beta_{123}M$ have positive scalar curvature
functions.

\medskip
{\bf \ref{s7}.2\qua Surgery theorem in the case of
manifolds without singularities}\qua Here we briefly review key results
on the connection between positive scalar curvature metric and surgery
for manifolds without singularities. The first basic result is due to
Gromov--Lawson \cite[Theorem A]{GL1} and to Schoen--Yau \cite{SY}. A
detailed ``textbook'' proof may be found in \cite[Theorem 3.1]{RS3}.
\begin{Theorem}[Gromov--Lawson \cite{GL1}, \ 
Schoen--Yau \cite{SY}]\label{rev:Th1} Let $M$ be a closed\break manifold,
not necessarily connected, with a Riemannian metric of positive scalar
curvature, and let $M^{\prime}$ is obtained from $M$ by a surgery of
codimension $\geq 3$. Then $M^{\prime}$ also admits a metric of
positive scalar curvature.
\end{Theorem}
To get started with $\Sigma$--manifolds we need an ``improved version''
of Theorem \ref{rev:Th1} which is due to Gajer \cite{Gaj}. 
\begin{Theorem}[Gajer \cite{Gaj}]\label{rev:Th2} 
Let $M$ be a closed manifold, not necessarily connected, with a
Riemannian metric $g$ of positive scalar curvature, and let
$M^{\prime}$ is obtained from $M$ by a surgery of codimension $\geq
3$. Then $M^{\prime}$ also admits a metric $g^{\prime}$ of positive
scalar curvature. Furthermore, let $W$ be the trace of this surgery
(ie, a cobordism $W$ with $\p W= M\sqcup -M^{\prime}$). Then there is
a positive scalar curvature metric $\bar{g}$ on $W$, so that
$\bar{g}=g+ dt^2$ near $M$ and $\bar{g}=g^{\prime}+ dt^2$ near
$M^{\prime}$. 
\end{Theorem}
In order to use the above Surgery Theorems, one has to specify certain
structure of manifolds under consideration. This structure (known as
$\gamma$--structure) is determined by the fundamental group $\pi_1(M)$,
and the Stiefel--Whitney classes $w_1(M)$, and $w_2(M)$. Indeed, it is
well-known that the fundamental group $\pi$ is crucially important for
the existence question. Then there is clear difference when a manifold
$M$ is oriented or not (which depends on $w_1(M)$). On the other hand,
a presence of the $Spin$--structure (which means that $w_2(M)=0$) gives
a way to use the Dirac operator on $M$ to control the scalar curvature
via the vanishing formulas. Stolz puts together those invariants to
define a $\gamma$--structure, see \cite{St3}.
In the case we are interested in, all manifolds are simply-connected
and $Spin$, thus we will state only a relevant Bordism Theorem (see,
say, \cite[Theorem 4.2]{RS3} for a general result).
\begin{Theorem}\label{rev:Th3}
Let $M$ be a simply connected $Spin$ manifold, $\dim M\geq 5$. Then
$M$ admits a metric of positive scalar curvature if and only if there
is some simply-connected $Spin$--manifold $M^{\prime}$ of positive
scalar curvature in the same $Spin$--bordism class.
\end{Theorem}
\vspace{2mm}

\noindent
{\bf \ref{s7}.3\qua Surgery theorem in the case of manifolds with
singularities}\qua 
Let $M$ be a $\Sigma$--manifold with $\Sigma=(P_i)$,
$(P_i,P_j)$ or $(P_i,P_j,P_k)$. Here $P_i$ are arbitrary closed
manifolds. Let $\dim M = n$, and $\dim P_i = p_i$, $i=1,2,3$. Then we
denote $\dim \beta_i M = n_{i} = n- p_i -1$, $\dim \beta_{ij} M =
n_{ij} = n - p_i -p_j -2$, and $\dim \beta_{123} M = n_{123} = n - p_1
-p_2 - p_3 -3$. The manifolds $\beta_i M$, $\beta_{ij} M $ and
$\beta_{ijk} M$ are called $\Sigma$--strata of $M$.

We say that a $\Sigma$--manifold $M$ is {\sl simply connected} if $M$
itself is simply connected and all $\Sigma$--strata of $M$ are simply
connected manifolds. 
\begin{Theorem}\label{rev:Th4}
Let $M$ be a simply connected $Spin$ $\Sigma$--manifold, $\dim M = n$,
so that all $\Sigma$--strata manifolds are nonempty, and satisfying the
following conditions:
\begin{enumerate}
\item[{\rm (1)}] if $\Sigma=(P_i)$, then $n -p_i \geq 6$; 
\item[{\rm (2)}] if $\Sigma=(P_i,P_j)$, then $n - p_i -p_j \geq 7$;  
\item[{\rm (3)}] if $\Sigma=(P_i,P_j,P_k)$, then $n - p_i -p_j -p_k \geq 8$.  
\end{enumerate}
Then $M$ admits a positive scalar curvature if and only if there is
some simply-connected $Spin^{\Sigma}$--manifold $M^{\prime}$ of
positive scalar curvature in the same $Spin^{\Sigma}$--bordism class.
\end{Theorem}
\begin{Remark}
The role of the manifolds $M$ and $M^{\prime}$ are not symmetric
here. For instance, it is important that $M$ has all $\Sigma$--strata
manifolds nonempty, however, the manifold $M^{\prime}$ may have empty
singularities.
\end{Remark}
\begin{proof}
(1)\qua Let $W$ be a $Spin^{\Sigma}$--cobordism between $M$ and
$M^{\prime}$.  Then $\beta_i W$ is a $Spin$--cobordism between $\beta_i
M$ and $\beta_i M^{\prime}$. By condition, $\beta_i M^{\prime}$ is
simply connected, and $\dim \beta_i M^{\prime} = \dim M \geq 5$.  We
notice that there is a sequence of surgeries on the manifold $\beta_i
W$ (relative to the boundary $\p\beta_i M^{\prime}$) so that the
resulting manifold is 2--connected (see an argument given in
\cite[Proof of Theorem A]{GL1}). Let $V$ be a trace of this
surgery. Then its boundary is decomposed as
$$
\p V = \beta_i W \cup (\beta_i M \times I) \cup (\beta_i M^{\prime} \times I)
\cup L_i.
$$
We glue together the manifolds $W$ and $- V\times P_i $:
$$
W^{\prime} := W \cup_{\beta_i W\times P_i} - V\times P_i .
$$
Then the boundary of $WW^{\prime}$ (as a $Spin^{\Sigma}$--manifold) is
$$
\begin{array}{l}
\strut0{10}\displaystyle
\delta W^{\prime} = \left(M\cup (\beta_i M \times I \times P_i )\right)
\sqcup
\left(M^{\prime}\cup (\beta_i M^{\prime} \times I \times P_i )\right)\cong
M\sqcup M^{\prime},
\\
\mbox{and} \ \ \ \beta_i W^{\prime} = L_i, \  \ \ 
\mbox{with} \ \ \ \p L_i = \beta_i M \sqcup \beta_i M^{\prime}.
\end{array}
$$
Now we use Theorem \ref{rev:Th2} to ``push'' a positive scalar
curvature metric from $\beta_i M^{\prime}$ through $L_i$ to $\beta_i
M$ keeping it a product metric near the boundary. At this point a
psc-metric $g_i$ on $L_i$ may be such that the product metric
$g_i\times g_{P_i}$ is not of positive scalar curvature. We find
$\epsilon>0$ so that the product metric $\epsilon g_i\times g_{P_i}$
has positive scalar curvature, and then we attach one more cylinder
$L_i\times P_i\times [0,a]$ with the metric 
$$
g_i(t):=\frac{a-t}{a} g_i\times g_{P_i} + \frac{t}{a}\epsilon  g_i\times g_{P_i} + dt^2.
$$
We use metric $g_i(t)$ to fit together the metric already constructed
on $W'^{\prime}$ with the metric on $L_i\times P_i\times [0,a]$.  In
particular, there is $a>0$ so that the restriction of $g_i(t)$ on
$\beta_i M^{\prime} \times P_i\times [0,a]$ has positive scalar
curvature (since an isotopy of positive scalar curvature metrics
implies concordance). By small perturbation, we can change $g_i(t)$, so
that it has positive scalar curvature and it is a product near the
boundary.  Then we do surgeries on the interior of $W^{\prime}$ to
make it 2--connected. Let $W^{\prime\prime}$ be the resulting
manifold. In particular, $\beta_i W^{\prime\prime}= \beta_i W^{\prime}
= L_i$.  Finally we use 
``push'' a positive scalar curvature metric from $M^{\prime}$ to $M$
through $W^{\prime\prime}$ keeping it a product metric near the
singular stratum $\beta_i \beta_i W^{\prime} = L_i$.

(2)\qua Let $M$ be a simply connected $Spin$ $\Sigma$--manifold, with
$\Sigma=(P_i,P_j)$, and $n - p_i -p_j \geq 7$. By condition, the
singular stratum $\beta_{ij} M\neq \emptyset$. Let $W$ be a
$Spin^{\Sigma}$--cobordism between $M$ and $M^{\prime}$. In particular,
we have $\p \beta_{ij} W = \beta_{ij} M \sqcup \beta_{ij}
M^{\prime}$. Recall that $\beta_{ij} W \times P_i\times P_j $ is
embedded to the union
$$
(\beta_{i} W \times P_i)\cup (\beta_{j} W \times P_j)
$$
together with the colors 
$$
\beta_{ij} W \times P_i\times P_j \times [-\epsilon, \epsilon]
$$ 
By conditions, the manifolds $\beta_{ij} M$ $\beta_{ij}
M^{\prime}$ are simply connected, and $\dim \beta_{ij} M$ $ = \dim
\beta_{ij} M^{\prime}\geq 5$. As above, there is a surgery on
$\beta_{ij} W$ (relative to the boundary $\p \beta_{ij} W = \beta_{ij}
M \sqcup \beta_{ij} M^{\prime}$) so that a resulting manifold is
2--connected. Let $V_{ij}$ be the trace of this surgery: 
$$
\p V =  \beta_{ij}  W \cup (\beta_{ij} M \times I) \cup (\beta_{ij} M^{\prime} \times I)
\cup L_{ij}.
$$
We glue together the manifolds 
$$
W \ \ \ \mbox{and} \ \ \ \ - V \times [-\epsilon, \epsilon]\times
P_i\times P_j
$$
to obtain a manifold $W^{\prime}$, where we identify
$$
\begin{array}{l}
\strut0{10}\ \beta_{ij} W \times P_i\times P_j \times [-\epsilon, \epsilon] \subset 
(\beta_{i} W \times P_i)\cup (\beta_{j} W \times P_j)\ \ \ \mbox{and}
\\
-\beta_{ij} W \times P_i\times P_j \times [-\epsilon, \epsilon] \subset 
-\p V \times [-\epsilon, \epsilon]\times
P_i\times P_j,
\end{array}
$$
see Figure 8. 

The resulting manifold $W^{\prime}$ (after smoothing corners and
extending metric according with the Surgery Theorem construction) is
such that $\beta_{ij} W^{\prime} = L_{ij}$ is 2--connected cobordism
between $\beta_{ij} M$ and $\beta_{ij} M^{\prime}$. Thus we can
``push'' a positive scalar curvature metric from
$\beta_{ij}M^{\prime}$ to $\beta_{ij} M$ through the cobordism
$\beta_{ij} W^{\prime}$. Thus we obtain a psc-metric $g_{ij}$ on
$\beta_{ij}M^{\prime}$ which is a product near boundary. In general,
the product metric $g_{ij}\times g_{P_i}\times g_{P_j}$ on $\beta_{ij}
W \times P_i\times P_j $ is not of positive scalar curvature. Then we have 
to attach one more cylinder 
$$
\beta_{ij} W^{\prime} \times [-\epsilon, \epsilon] \times I \times P_i\times P_j 
$$
to ``scale'' the metric $g_{ij}\times g_{P_i}\times g_{P_j}$ to a
positive scalar curvature metric $\epsilon_{ij} g_{ij}\times
g_{P_i}\times g_{P_j}$ through an appropriate homotopy. 

\hspace*{15mm}\PSboxa <0pt,10pt> {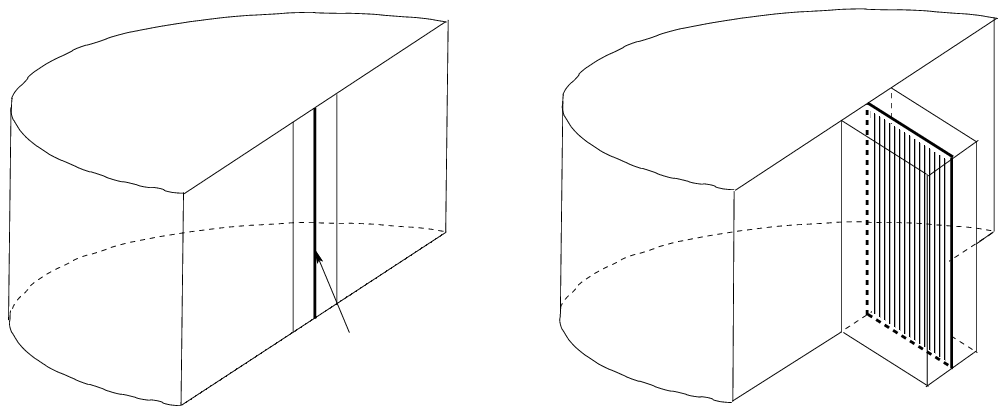}{30mm}{45mm}
\begin{picture}(10,10)
\put(-85,55){{\small $W$}}
\put(-65,25){{\small $M$}}
\put(-65,85){{\small $M^{\prime}$}}
\put(-10,20){{\small $\beta_{ij} W \times P_i\times P_j $}}
\put(130,5){{\small $-V \times [-\epsilon, \epsilon]\times
P_i\times P_j$}}
\end{picture}

\centerline{{\small Figure 8}}

Then we consider the manifolds $\beta_i W^{\prime}$ and $\beta_j
W^{\prime}$.  Again, we perform surgeries on the interior of $\beta_i
W^{\prime}$, $\beta_j W^{\prime}$ to get 2--connected manifolds $L_i$
and $L_j$. Let $V_i$, $V_j$ be the traces of these surgeries:
$$
\p V_i = \beta_i W^{\prime} \cup \beta_{ij} W^{\prime} \times P_j \cup L_i, 
\ \ \ 
\p V_j = \beta_i W^{\prime} \cup \beta_{ij} W^{\prime} \times P_i \cup L_j.
$$
Now we attach the manifolds $-V_i\times P_i$ and $-V_j\times V_j$ to
$W^{\prime}$ by identifying
$$
\begin{array}{lll}\strut0{10}
\beta_i W^{\prime} \times P_i \subset W^{\prime} \ \ &\mbox{and} \ \  \
\beta_j W^{\prime} \times P_j \subset W^{\prime} \ \ \ &\mbox{with}
\\
-\beta_i W^{\prime} \times P_i \subset -\p V_i, \ \ \ &\mbox{and} \ \ \ 
-\beta_j W^{\prime} \times P_j \subset -\p V_j &
\end{array}
$$
respectively. Let $W^{\prime\prime}$ be the resulting manifold (after
an appropriate smoothing and extending a metric), see Figure 9. Notice
that $W^{\prime\prime}$ is still a $Spin^{\Sigma}$--cobordism between
$M$ and $M^{\prime}$.

This procedure combined with an appropriate metric homotopy gives
$W^{\prime\prime}$ together with a metric $g^{\prime\prime}$ on
$W^{\prime\prime}$, so that it is a product metric near the boundary,
its restriction on $M^{\prime}$ has positive scalar curvature, and its
restriction on the manifolds
$$
\beta_{i} W^{\prime\prime}\times P_i, \ \ \beta_{j} W^{\prime\prime}\times P_j, \
\ \beta_{i} W^{\prime\prime}\times [-\epsilon,\epsilon]\times P_i\times P_j
$$
are psc-metrics
$$
g_i\times g_{P_i}, \  \ g_j\times g_{P_j}\  \ g_{ij} \times g_{P_i}\times g_{P_j} +dt^2
$$
respectively (for some psc-metrics $g_i$, $g_j$, $ g_{ij}$). It
remains to perform surgeries on the interior of $ W^{\prime\prime}$ to
get a 2--connected manifold, and finally push a psc-metric from
$M^{\prime}$ to $M$ relative to the boundary. 

\hspace*{15mm}\PSbox{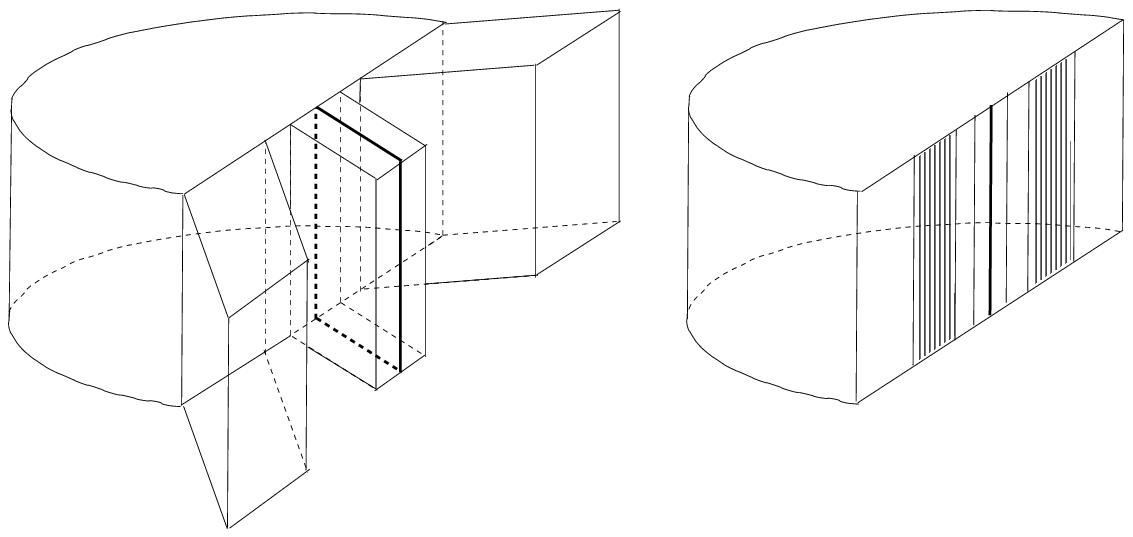}{30mm}{53mm}
\begin{picture}(0,10)
\put(-65,45){{\small $M$}}
\put(-65,105){{\small $M^{\prime}$}}
\put(-90,-15){{\small $W^{\prime}\cup (-V_i\times P_i)\cup (-V_j\times P_j)$}}
\put(0,15){{\small $-V_i\times P_i$}}
\put(45,60){{\small $-V_j\times P_j$}}
\put(150,15){{\small $W^{\prime\prime}$}}
\end{picture}
\vglue 20pt

\centerline{{\small Figure 9}}

A proof of (3) is similar. 
\end{proof}

\section{Proof of Theorem \ref{ThA}}\label{s8}
First we recall the main construction from \cite{St1}.  Let $G=PSp(3)
= Sp(3)/(\Z/2)$, where $\Z/2$ is the center of $Sp(3)$.  Let $g_0$ be
the standard metric on $\H\P^2$. Recall that the group $G$ acts on
$\H\P^2$ by isometries of the metric $g_0$. Let $E \longrightarrow  B$ be a
geometric $\H\P^2$--bundle, ie, $E \longrightarrow  B$ is a bundle with a fiber
$\H\P^2$ and structure group $G$. Each geometric $\H\P^2$--bundle $E
\longrightarrow  B$ is given by a map $f \co B \longrightarrow  BG$ by taking first the
associated principal $G$--bundle, and then by ``inserting'' $\H\P^2$ as
a fiber employing the action of $G$.  Assume that $B$ is a $Spin$
manifold. Then the correspondence $(B,f)\mapsto E$ gives the
homomorphism $T \co \Omega_{n-8}^{Spin}(BG) \longrightarrow  \Omega_{n}^{Spin}$.
Let $X$ be a finite $CW$--complex.  The homomorphism $T$ actually gives
the transformation
$$
T \co \Omega_{n-8}^{Spin}(X\wedge BG_+) \longrightarrow  \Omega_{n}^{Spin}(X),
$$
which may be interpreted as the transfer, and induces the map at the
level of the classifying spectra
\begin{equation}\label{pr1}
T \co MSpin\wedge \Sigma^8BG_+ \longrightarrow  MSpin,
\end{equation}
see details in \cite{St1}. Consider the composition
\begin{equation}\label{pr2}
\begin{CD}
\Omega_{n-8}^{Spin}(X\wedge BG_+) @>T>> \Omega_{n}^{Spin}(X)
@>\alpha>> ko_n(X)
\end{CD}
\end{equation}
Here is the result due to S Stolz, \cite{St2}:
\begin{Theorem}\label{Stolz} Let $X$ be a $CW$--complex. Then
there is an isomorphism at the 2--local category: {\rm 
$$
ko_n(X)\cong \Omega_{n}^{Spin}(X) / \Im T 
$$
}
\end{Theorem}
Let $\Sigma=\Sigma_1$ or $\Sigma_2$, or $\Sigma_3$ or $\eta$, and
$X_{\Sigma}$ be the corresponding spectrum, so that
$MSpin^{\Sigma}\cong MSpin\wedge X_{\Sigma}$. The map 
$T \co MSpin\wedge \Sigma^8BG_+ \longrightarrow  MSpin$ induces the
map 
$$
T^{\Sigma} \co MSpin\wedge \Sigma^8 BG_+ \wedge X_{\Sigma} \longrightarrow  MSpin
\wedge X_{\Sigma}.
$$
Consider the composition
$$
\begin{CD}
MSpin\wedge \Sigma^8 BG_+ \wedge X_{\Sigma} @>T^{\Sigma} >> MSpin
\wedge X_{\Sigma} @>\alpha^{\Sigma} >> ko\wedge X_{\Sigma}.
\end{CD}
$$
We use Theorems \ref{Sigma1} and \ref{Sigma3} to derive the following
conclusion from Theorem \ref{Stolz}.
\begin{Corollary}\label{pr3}
Let $\Sigma=\Sigma_1$ or $\Sigma_2$, or $\Sigma_3$ or $\eta$. 
Then
there is an isomorphism at the 2--local category:
$$
ko^{\Sigma}_n \cong MSpin^{\Sigma}_n/\Im T^{\Sigma}.
$$
\end{Corollary}
We remind here that the homomorphism $ko^{\Sigma}_n \to KO^{\Sigma}_n$
is a monomorphism for $n\geq 0$.  

Corollary \ref{pr3} describes the situation in $2$--local category.
Now we consider what is happening when we invert $2$. Consider first
the case when $\Sigma=\Sigma_1$.  Then we have a cofibration:
$$
\begin{CD}
MSpin [\frac{1}{2}] @>\cdot 2 >>  MSpin[\frac{1}{2}] @>>> MSpin^{\Sigma_1}[\frac{1}{2}]
\end{CD}
$$
Clearly the map $\cdot 2 \co \ MSpin[\frac{1}{2}] \to MSpin[\frac{1}{2}] $
is a homotopy equivalence. Thus $MSpin^{\Sigma_1}[\frac{1}{2}]\cong
pt$. The case $\Sigma=(\eta)$ is more interesting. Here we have the
cofibration:
$$
\begin{CD}
S^1 \wedge MSpin[\frac{1}{2}] @>\cdot\eta>> MSpin[\frac{1}{2}] @>\pi^{\eta}>> 
MSpin^{(\eta)}[\frac{1}{2}].
\end{CD}
$$
Notice that $\eta=0$ in $p$--local homotopy $\Omega_1^{Spin}\otimes
\Z[\frac{1}{2}]$. Thus we have a short exact sequence:
$$
\begin{CD}
0\to \Omega_n^{Spin}\otimes\Z[\frac{1}{2}] @>\pi^{\eta}_*>> 
\Omega_n^{Spin,\eta}\otimes\Z[\frac{1}{2}]  @>\beta>>
\Omega_{n-2}^{Spin}\otimes\Z[\frac{1}{2}]  \to 0.
\end{CD}
$$
This sequence has very simple geometric interpretation. Let $w\in
\Omega_2^{Spin,\eta}\cong \Z$ be an element represented by an
$(\eta)$--manifold $W$, so that $\beta W = 2$.  Let $t= \frac{w}{2}\in
\Omega_2^{Spin,\eta} \otimes\Z[\frac{1}{2}]$. 

Now let $c\in \Omega_n^{Spin,\eta}\otimes\Z[\frac{1}{2}]$, and $\beta
c = b$. Let $a= c - tb$, then $c = a +tb$ for $a\in 
\Omega_n^{Spin}\otimes\Z[\frac{1}{2}]$, $b \in
\Omega_{n-2}^{Spin}\otimes\Z[\frac{1}{2}]$ for any element $c\in
\Omega_n^{Spin,\eta }\otimes\Z[\frac{1}{2}]$. Furthermore, this
decomposition is unique once we choose an element $t$.  Recall that
$\Omega_*^{Spin}\otimes\Z[\frac{1}{2}] \cong
\Omega_*^{SO}\otimes\Z[\frac{1}{2}]$ is a polynomial algebra $
\Z[\frac{1}{2}][x_1,x_2,\ldots, x_j, \ldots]$ with $\deg x_j =4j$.
Thus we obtain
\begin{equation}\label{rev:eq1}
\begin{array}{l}
\Omega_n^{Spin,(\eta)}\otimes\Z[\frac{1}{2}] \cong
\{
\begin{array}{ll}
\Omega_n^{Spin}\otimes\Z[\frac{1}{2}] & \mbox{if} \ n=4k 
\\
\strut{15}5\Omega_{n-2}^{Spin}\otimes\Z[\frac{1}{2}] & \mbox{if} \ n=4k+2
\\
0  & \mbox{otherwise}
\end{array}
\right.
\end{array}
\end{equation}
Furthermore, as it is shown in \cite[Proposition 4.2]{KS} there are
generators $x_j=[M^{4j}]$ of the polynomial algebra
$\Omega_*^{Spin}\otimes\Z[\frac{1}{2}]$, so that the manifolds
$M^{4j}$ are total spaces of geometric $\H\P^2$--bundles (for all $j\geq
2$). In particular, it means that the groups
$\Omega_{4j}^{Spin}\otimes\Z[\frac{1}{2}]$ are in the ideal $\Im T
\subset \Omega_{*}^{Spin}\otimes\Z[\frac{1}{2}]$. Now the formula
(\ref{rev:eq1}) shows that the groups
$\Omega_n^{Spin,(\eta)}\otimes\Z[\frac{1}{2}] $ are in the ideal $\Im
T^{\eta}$.  We obtain the isomorphism in integral homotopy groups: $
ko^{\Sigma}_n \cong \Omega^{Spin,\eta}_n/\Im T^{\eta} $.

The cases $\Sigma =\Sigma_2$ or $\Sigma =\Sigma_3$ are similar to the
case $\Sigma =\Sigma_1$: here we have that
$MSpin^{\Sigma_i}[\frac{1}{2}]\cong pt$ for $i=2,3$.
 
Thus in all cases we conclude that any element $x\in \Im T^{\Sigma}$
may be represented by a simply connected $\Sigma$--manifold admitting a
psc-metric. Here the restriction that $\dim x \geq d(\Sigma)$ is
essential. Thus we conclude that if a simply-connected $Spin$ $\Sigma$
manifold $M$ with $\dim M \geq d(\Sigma)$ is such that $[M]\in \Ker
\alpha^{\Sigma}$, then $M$ admits a psc-metric.

Now we prove the necessity.

Let $M$ be a simply connected $Spin$--manifold of dimension
$\dim M \geq d(\Sigma)$. What we really must show is that if there is
a psc-metric on $M$, then $\alpha^{\Sigma}([M])=0$ on the group
$KO^{\Sigma}_*$.

The case $\Sigma_1 = (P_1)=\<2\>$ is done in \cite{R}, where it is
shown that $\alpha^{\Sigma}([M]) \in KO^{\<2\>}_*$ coincides with the
index of the Dirac operator on $M$, and that the index
$\alpha^{\Sigma}([M])$ vanishes if $M$ has a psc-metric.

The next case to consider is when $\Sigma = \eta=(P_2)$. Let $M$ be a
closed $\Sigma$--manifold, ie, $\p M = \beta_2 M \times P_2$, where
$P_2$ is a circle with the nontrivial $Spin$ structure. Let $g$ be a
psc-metric on $M$. In particular, we have a psc-metric $g_{\beta_2
M}$.  Then, as we noticed earlier, the circle $P_2$ is zero-cobordant 
as $Spin^c$--manifold. More precisely, we choose a disk $D^2$ with
$\p D^2 = P_2$, and construct the manifold
$$
\ov{M} = M\cup\(-\beta_2 M \times D^2\)
$$
where we identify $\p M \beta_2 M \times P_2$ with $\p(-\beta_2 M
\times D^2)$. There is a canonical map 
$$
h \co \ov{M} \longrightarrow  \C\P^1
$$
which sends $M\subset \ov{M}= M\cup\(-\beta_2 M \times D^2\)$ to the point,
and $\beta_2 M \times D^2$ to $\C\P^1=S^2$ by the composition
$$
\beta_2 M \times D^2 \longrightarrow  D^2 \longrightarrow  D^2/S^1 =\C\P^1. 
$$
The map $h$ composed with the inclusion $\C\P^1 \subset \C\P^{\infty}$
gives the map $\bar{h} \co \ov{M} \longrightarrow \C\P^{\infty}$, and,
consequently, a linear complex bundle $\xi \longrightarrow \ov{M}$
which is trivialized over $M$. The $Spin$--structure on $M$ together
with the linear bundle $\xi \longrightarrow M$ determines a
$Spin^c$--structure on $\ov{M}$. To choose a metric $g_0$ on the disk,
We identify $D^2$ with the standard hemisphere $S^2_+$ with a small
color attached to the circle $S^1$, so that the metric $g_0|_{S^1}$ is
the standard flat metric $d\theta^2$. Then we have the product metric
$g_{\beta_2 M}\times g_0$ on $\beta_2 M \times D^2$.
Together with the metric $g$ on $M$, it gives a
psc-metric $\ov{g}$ on $\ov{M}$. We choose a $U(1)$--connection on the
linear bundle $\xi \longrightarrow M$, and let $F$ be its curvature
form. We notice that since $\xi$ is trivialized over $M\subset
\ov{M}$, the form $F$ is supported only on the submanifold $-\beta_2 M
\times D^2\subset \ov{M}$. Moreover, we have defined the bundle $\xi
\longrightarrow \ov{M}$ as a pull-back from the tautological complex
linear bundle over $\C\P^1$. Thus locally we can choose a basis
$e_1,e_2,\ldots e_n$ of the Clifford algebra, so that $F(e_1,e_2)\neq
0$, and $F(e_i ,e_j) =0$ for all other indices $i,j$. Notice also that
the scalar curvature function $R_{\ov{g}}= R_{g_{\beta_2 M}} + R_{g_0}$.
Let $D$ be the Dirac operator on the canonical bundle $S(\ov{M})$ of
Clifford modules over $\ov{M}$. We have the BLW--formula
\begin{equation}\label{pr4}
D^2 = \nabla^*\nabla +{1\over 4} \left(R_{\ov{g}} + R_{g_0}\right) +
{1\over 2} F(e_1,e_2)\cdot e_1 \cdot e_2.
\end{equation}
Now we scale the metric $g_{\beta_2 M}$ to the metric $\epsilon^2 \cdot g_{\beta_2 M}$
with the scalar curvature $R_{\epsilon^2 \cdot g_{\beta_2 M}} =
\epsilon^{-2} \cdot R_{g_{\beta_2 M}}$. Clearly this scaling does not
effect the connection form since the scaling is in the ``perpendicular
direction''. Let $\epsilon>0$ be such that the term 
$$
{1\over 4} \left(\epsilon^{-2} \cdot R_{g_{\beta_2 M}} + R_{g_0}\right)
$$ 
will dominates the connection term ${1\over 2} F(e_1,e_2)\cdot e_1
\cdot e_2$. Then we attach the cylinder $\beta_2 M \times [0,a] \times D^2$ (for
some $a>0$) with the metric $g_{\beta_2 M}(t) \times g_0$, where  
$$
g_{\beta_2 M}(t) 
= \frac{a-t}{a} g_{\beta_2 M}  + 
\frac{t}{a}\epsilon^2 _{\beta_2 M}  + dt^2,
$$
so that the metric $g_{\beta_2 M}(t) \times g_0$ has positive scalar
curvature, and is a product metric near the boundary.  Thus with that
choice of metric, the right-hand side in (\ref{pr4}) becomes positive,
which implies that the Dirac operator $D$ is invertible, and hence
$\ind (D)\in K_*$ vanishes. This completes the case of
$\eta$--singularity.
\begin{Remark}
Here the author would like to thank S Stolz for
explaining this matter.
\end{Remark}
The case $\Sigma_2=(P_1,P_2)$ is just a combination of the above
argument and the BLW--formula for $Spin^c$ $\Z/k$--manifolds given by Freed
\cite{Freed}. 

The last case, when $\Sigma=\Sigma_3=(P_1,P_2,P_3)$ there is nothing to prove
since $KO^{\Sigma_3}$ is a contractible spectrum, and thus any
$\Sigma_3$--manifold has a psc-metric. Indeed, we have that
$$
\Omega^{Spin,\Sigma_3}_n \cong \Im T ,  \ \ \ \mbox{if} \ \ \ n\geq 17.
$$
This completes the proof of Theorem \ref{ThA}.

\end{document}